\documentclass[11pt]{amsart}
\usepackage{amsmath}
\usepackage{amssymb}
\usepackage[utf8]{inputenc}
\usepackage[T1]{fontenc}
\usepackage{indentfirst}
\usepackage{amsthm}
\usepackage{mathrsfs}
\usepackage{moreverb}
\usepackage{dsfont}
\usepackage{amsfonts}
\usepackage{xcolor}
\usepackage{url}
\usepackage{mathtools}
\usepackage{caption}
\usepackage{setspace}
\usepackage{subcaption}
\usepackage{float}
\usepackage{bigints}

\newtheorem{remark}{Remark}
\newtheorem{notation}{Notation}
\newtheorem{assumption}{Assumption}
\newsavebox\myboxA
\newsavebox\myboxB
\newlength\mylenA
\newcommand{\keyword}{\textit{Keywords.}}

\begin{document}
\title[]{Efficient Asymptotic Models for Axisymmetric Eddy Current Problems in Linear Ferromagnetic Materials}
\author[D. ABOU EL NASSER EL YAFI and V. PÉRON]{Dima ABOU EL NASSER EL YAFI, Victor PÉRON}

\begin{abstract}
    The problem under consideration is that of time-harmonic eddy current problems in linear ferromagnetic materials surrounded by a dielectric medium with a smooth common interface. Assuming axisymmetric geometries and orthoradial axisymmetric data, we construct an efficient multiscale expansion for the orthoradial solution that provides reduced computational costs. We investigate numerically the accuracy of the approach using an analytical procedure and infinite cylinders as well. It results that the computation of two asymptotics is sufficient to ensure accurate solutions in the case of low frequencies.
    \\ \noindent \keyword \
    Multiscale expansion; Eddy current problems; Analytical method; Ferromagnetic materials; Axisymmetric geometry 
    \end{abstract}

\maketitle

\section{Introduction}
 Eddy currents arise due to the time varying magnetic field crossing metals \cite{biro1999edge}. The distribution of the current density in this case is restricted at a boundary layer near the metallic surface, and diminishes exponentially inside the conducting medium. This phenomenon is called the \textit{skin effect} \cite{rytov1940calcul,stephan1983solution,maccamy1987skin,bossavit2004electromagnetisme,caloz2011influence}. Eddy currents generate energy losses that have two sided effect in the industrial field \cite{rodriguez2010eddy}. On the one hand, these currents can have a good use such as induction heating or for the design of electromagnetic breaking systems. On the other hand, eddy currents can also produce "undesirable" power losses in the form of heating for example which can affect the performance of some electrical devices. Summing up, studying eddy currents is crucial for engineering applications in electromagnetism.

The mathematical and numerical analysis of the eddy current problems have been the interest of many works during the past decades \cite{maccamy1984solution,krahenbuhl1993thin,biro1999edge,buffa2000justification,hiptmair2002symmetric,costabel2003singularities,bossavit2004electromagnetisme,bermudez2007transient,rodriguez2010eddy}. Because of the small skin depth inside the conductors, the classical numerical methods are challenging to apply. To overcome this difficulty, it is possible to develop an asymptotic method that derives approximate models with less computational costs. The asymptotic approach is often employed for physical problems involving a small or large parameter. This method gives an accurate approximation of the problem by solving an ordered sequence of subproblems independent of the latter parameter.

We refer the reader to \cite{rytov1940calcul,hariharan1982integral,maccamy1984solution,haddar2008generalized,dauge2010comportement,caloz2011influence,perrussel2013asymptotic,schmidt2013unified,peron2019asymptotic,issa2019boundary} for previous works devoting to the asymptotic procedure in electromagnetic problems. For example, in \cite{hariharan1982integral,maccamy1984solution} authors investigated eddy current problems in the case of metals having infinite conductivity by applying first a boundary integral procedure and then an asymptotic procedure that reflects the skin effect in metals in both bi-dimensional and three-dimensional domains. Moreover, recent studies \cite{haddar2008generalized,dauge2010comportement,caloz2011influence} analyzed theoretically and numerically the electromagnetic field solution for the Maxwell equations through an asymptotic expansion for large conductivities. It is worthwhile to note that these previous works tackle the equations of electromagnetism set on a domain made of a dielectric and a non-magnetic conducting subdomains with a smooth common interface. On the other hand, several works have investigated the asymptotic approach for eddy current problems in a bi-dimensional setting where the conducting medium is non-magnetic and has a corner singularity on the conductor-dielectric interface \cite{buret2012eddy,dauge2014corner,dular2014perfect}. These works have shown that the asymptotic approach was strongly affected by adding corrections, especially near corners, in order to obtain accurate asymptotic models.

This paper continues a study begun in \cite{peron2021magnetic} and \cite{abounumerical} of the time harmonic eddy current problems in linear ferromagnetic materials with a smooth interface. The work in \cite{peron2021magnetic} was restricted for the theoretical results whereas in \cite{abounumerical} was concerned essentially with numerical validation of the asymptotic procedure employed in \cite{peron2021magnetic}. In both cases, the study was restricted to a very special class of two-dimensional problems using a multiscale expansion. Besides, \cite{abounumerical} is not a straightforward application of \cite{peron2021magnetic}. More precisely, we identified in \cite{abounumerical} efficient asymptotic models, slightly different than those established in \cite{peron2021magnetic}, that provide reduced computational costs in time and memory allocation for a wide range of physical parameters. This present work treats a three-dimensional situation in axisymmetric geometry.

However, three dimensional computations can be very expensive. In a number of cases, it is possible to reduce the problem by assuming that the geometry is invariant by translation or rotation \cite{bernardi1999spectral,assous2003solution,caloz2011influence,ciarlet2011numerical}. In this context, we choose to consider a special class of axisymmetric geometry that reduces our problem to a one-dimensional scalar model. Our analysis is twofold. First, we present elements of derivation for the multiscale expansion of the one-dimensional solution near the conductor-insulator interface. We identify efficient asymptotic models that reduce the computational costs. Then, we evaluate the performance of the resulting models by presenting numerical results.

In this work we assess the performance of our efficient asymptotic models analytically in the case of infinite cylinders. Indeed, our analytical procedure follows the spirit of \cite{bermudez2007transient} where authors tackled eddy current problems for large conductivities and for the case of infinite cylinders as well. There are many differences between our paper and \cite{bermudez2007transient}. For example, in \cite{bermudez2007transient} authors analyze the performance of the finite element method (FEM) applied to the considered eddy current problems. However, our numerical study is based on analytical methods in order to highlight the good accuracy of our asymptotic approach, for an example of unbounded domain. Actually, the FEM was previously applied for the eddy current problems in linear ferromagnetic materials \cite{abounumerical} where their analytical solution was not obvious to calculate because of the complexity of the considered bounded geometry. Moreover, in \cite{bermudez2007transient}, authors considered infinite cylinders, in width and length, consisting of a core material surrounded by a crucible and an extremely thin coil. The crucible itself is made of several concentric layers with different materials. In our case, for the sake of simplicity, we consider only two different layers: a ferromagnetic material surrounded by a dielectric material with a common smooth interface. Finally, we perform numerically a comparison of our asymptotic approach with the impedance method \cite{leontovich1948approximate,haddar2008generalized}.

The presentation of the paper proceeds as follows. Section \ref{section2} introduces the framework as well as the boundary value problem. In section \ref{section4}, we restrict our work to axisymmetric domains and orthoradial axisymmetric data. Moreover, we apply in the latter section a multiscale expansion for the orthoradial component of the magnetic vector potential and we identify efficient asymptotic models up to the order two. In section \ref{section6}, we present numerical results to assess the performance of the proposed models. Concluding remarks and perspectives are given in section \ref{section7}. In appendix \ref{Appendix A}, we provide elements of proof of the multiscale expansion given in the subsection \ref{subsection Multiscale expansion}. Appendix \ref{Appendix B} is dedicated to a deep calculation of the analytical solutions introduced in the section \ref{section 5}.  

\section{Problem setting} \label{section2}
\begin{figure}[h]
    \centering
    \includegraphics[width=0.5 \textwidth]{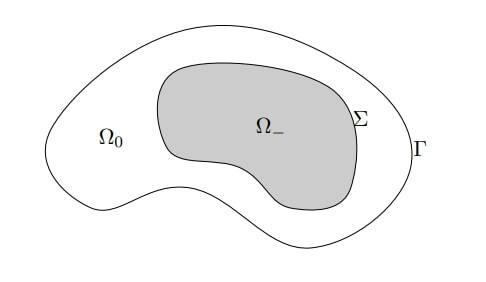}
    \caption{A cross section of the domain $\Omega$ and its subdomains $\Omega_{-}$, $\Omega_{0}$}
    \label{Figure 1}
\end{figure}

Throughout  the paper we denote by $\Omega \subset \mathbb{R}^{3}$ a smooth and connected domain with boundary $\Gamma,$ and $\Omega_{-}$ a smooth connected subdomain of $\Omega$ with boundary $\Sigma$. We denote by $\Omega_{0}$ the complementary of $\overline{\Omega}_{-},$ see for instance Figure \ref{Figure 1}.

\subsection{Notations and physical parameters}
We suppose that $\Omega_{0}$ is a dielectric medium which we consider for the sake of simplicity the free space, and $\Omega_{-}$ is a ferromagnetic material. The magnetic permeability and the conductivity are given by the following piecewise-constant functions $\underline{\mu}$ and $\underline{\sigma}$ respectively:
\begin{equation}
    \begin{array}{lllll}
     \underline{\mu}= \left\{\begin{array}{llll}
          \mu_{0} & & \mathrm{in}  & \Omega_{0}  \\
           \mu_{r} \mu_{0} & & \mathrm{in} & \Omega_{-}
      \end{array}\right.    
      & & \mathrm{and}  & & \underline{\sigma}= \left\{\begin{array}{llll}
          0 & & \mathrm{in}  & \Omega_{0}  \\
           \sigma > 0 & & \mathrm{in}  & \Omega_{-},
      \end{array}\right.
      \end{array}
\end{equation}
where $\mu_{0}=4 \pi \times 10^{-7}$ [H/m](henry per meter) and the relative permeability $\mu_{r}$ is assumed to be a large parameter. The angular frequency is denoted by $\omega > 0.$ In our work, $\omega$ and $\sigma$ are given parameters. We denote by $J_{s}$ the current source which is supposed for the sake of simplicity divergence free that is $\mathrm{div} \ J_{s} = 0 \ in \ \Omega,$ smooth enough and the support of $J_{s}$ does not meet $\Omega_{-}.$ 
We consider the following notations.
\begin{notation}
We denote by $h^{+}$ (resp. $h^{-}$) the restriction of any function $h$ in $\Omega_{0}$ (resp. $\Omega_{-}$).
\end{notation}

\begin{notation}
In order to introduce our asymptotic method, we define a small parameter $\varepsilon$ as follows
$$ \varepsilon=\dfrac{1}{\mu_{r} \delta}, $$
where $\delta$ is the skin depth and given by $$\delta = \sqrt{\dfrac{2}{\omega \sigma \mu_{r} \mu_{0}}}.$$
\end{notation}
\subsection{Boundary value problem}
The magnetic vector potential $\mathcal{A}=(\mathcal{A}^{+}, \mathcal{A}^{-})$ satisfies the following boundary value problem \cite{dular1994modelisation}
\begin{small}
	 \begin{equation} \label{mag pot A}
	 \begin{array}{ll}
	 \left\{ 
	 \begin{array}{ll}
	 \mathrm{curl \ curl \ \mathcal{A}^{+}}=\mu_{0} J_{s} &  \mathrm{in} \ \Omega_{0},\\ \\
	 \rm curl \ curl \ \mathcal{A}^{-} - i \ \omega \sigma \mu_{0} \mu_{r} \mathcal{A}^{-} = 0  &  \mathrm{in} \  \Omega_{-},\\ \\
	 \rm div \ \mathcal{A}^{-} = 0 &  \mathrm{in} \  \Omega_{-}, \\ \\
	 \rm \mathcal{A}^{+} \times \textit{n} = \mathcal{A}^{-} \times \textit{n}  &  \mathrm{on} \ \Sigma, \ \ \ \iff \\ \\
	 \rm curl \ \mathcal{A}^{+} \times \textit{n} = \mu_{r}^{-1} (curl \ \mathcal{A}^{-} \times \textit{n}) &  \mathrm{on} \ \Sigma,\\ \\
	 \rm \mathcal{A}^{-} \cdot \textit{n} = 0  &  \mathrm{on} \ \Sigma,\\ \\
	 \rm \mathcal{A}^{+} \times \textit{n} = 0   &  \mathrm{on} \ \Gamma. \\ \\
	 \end{array}
	 \right. & \left\{ \begin{array}{ll}
	 \mathrm{curl \ curl} \ \mathcal{A}^{+}=\mu_{0} J_{s} &  \mathrm{in} \  \Omega_{0},\\ \\
	 \rm curl \ curl \ \mathcal{A}^{-} - 2i\delta^{-2} \mathcal{A}^{-} = 0  &  \mathrm{in} \  \Omega_{-},\\ \\
	 \rm div \ \mathcal{A}^{-} = 0 &  \mathrm{in} \ \Omega_{-},\\ \\
	 \rm \mathcal{A}^{+} \times \textit{n} = \mathcal{A}^{-} \times \textit{n}  &  \mathrm{on} \ \Sigma, \\ \\
	 \rm curl \ \mathcal{A}^{+} \times \textit{ n} = \delta \varepsilon (curl \ \mathcal{A}^{-} \times \textit{n}) & \mathrm{on} \ \Sigma,\\ \\
	 \rm \mathcal{A}^{-} \cdot \textit{n} = 0  &  \mathrm{on} \ \Sigma,\\ \\
	 \rm \mathcal{A}^{+} \times \textit{n} =0   &  \mathrm{on} \  \Gamma. \\ \\
	 \end{array}
	 \right.
	 \end{array}
	 \end{equation}
\end{small}
This problem has to be completed by the gauge conditions \cite{buffa2000justification,hiptmair2002symmetric,costabel2003singularities,peron2019asymptotic}
\begin{equation} \label{gauge conditions}
   \begin{array}{lll}
     \mathrm{div} \ \mathcal{A}^{+} = 0 \ \mathrm{in} \ \Omega_{0} & \mathrm{and}  &   \int_{\Sigma} \mathcal{A}^{+} \cdot \textit{n} \ \mathrm{d}S=0.\\
   \end{array} 
\end{equation}

For numerical purposes, we used here the modified magnetic vector potential \cite[section 4.4 - page 25]{dular1994modelisation}, so its boundary value problem (\ref{mag pot A}) is deduced from that of the electric field.

\subsection{Variational formulation}
We define the following space 
\begin{equation}
    \mathbf{H}_{0}(\mathrm{curl}, \Omega)=\{ \mathrm{u} \in \mathbf{L}^{2}(\Omega) \ | \  \mathrm{curl} \ \mathrm{u} \in \mathbf{L}^{2}(\Omega), \ \mathrm{u} \times n=0 \ \mathrm{on} \ \Gamma \}.
\end{equation}
The variational space is the Hilbert space \textbf{Y}:
\begin{equation}
    \mathbf{Y}=\{ \mathrm{u} \in \mathbf{H}_{0}(\mathrm{curl}, \Omega) \ | \ \mathrm{div}\ \mathrm{u}^{+} \in \mathrm{L}^{2}(\Omega_{0}), \ \mathrm{div} \ \mathrm{u}^{-} \in \mathrm{L}^{2}(\Omega_{-}),\ \int_{\Sigma} \mathrm{u}^{+} \cdot n \ \mathrm{dS}=0 \}
\end{equation}
endowed with the norm 
\begin{equation*}
    \lVert \mathrm{u} \rVert_{\mathbf{Y}}^{2}=\lVert \mathrm{u} \rVert^{2}_{0,\Omega}+\lVert \mathrm{curl} \ \mathrm{u} \rVert^{2}_{0,\Omega}+\lVert \mathrm{div} \ \mathrm{u}^{+} \rVert_{0,\Omega_{0}}^{2}+\lVert \mathrm{div} \ \mathrm{u}^{-} \rVert_{0,\Omega_{-}}^{2}.
\end{equation*}
We introduce the small parameter $\nu=\dfrac{1}{\sqrt{\mu_{r}}}$ in the problem below.
For all $\nu > 0,$ the variational problem writes \\
\textit{Find} $\mathcal{A} \in \mathbf{Y}$ \textit{such that for all} $\mathrm{v} \in \mathbf{Y},$
\begin{equation} \label{VF A}
    a_{R}(\mathcal{A}, \mathrm{v})=\mu_{0}\int_{\Omega} J_{s} \cdot \overline{\mathrm{v}} \ \mathrm{dx}.
\end{equation}
Here the sesquilinear form in its regularized version $a_{R}$ is defined as 
\cite{costabel2003singularities,peron2019asymptotic}
\begin{equation*}
    \begin{array}{lll}
     a_{R}(\mathrm{u}, \mathrm{v}) & = & \nu^{2} \int_{\Omega_{-}} \mathrm{curl} \ \mathrm{u}^{-} \cdot \mathrm{curl} \ \overline{\mathrm{v}^{-}} \ \mathrm{dx} + \int_{\Omega_{0}} \mathrm{curl} \ \mathrm{u}^{+} \cdot \mathrm{curl} \ \overline{\mathrm{v}^{+}} \  \mathrm{dx}  + \int_{\Omega_{0}}  \mathrm{div} \ \mathrm{u}^{+} \  \mathrm{div} \ \overline{\mathrm{v}^{+}} \ \mathrm{dx} \\ & & + \int_{\Omega_{-}}  \mathrm{div} \ \mathrm{u}^{-} \  \mathrm{div} \ \overline{\mathrm{v}^{-}} \ \mathrm{dx} - \mathrm{i} \omega \sigma \mu_{0} \int_{\Omega_{-}} \mathrm{u}^{-} \cdot \overline{\mathrm{v}^{-}} \  \mathrm{dx}.
    \end{array} 
    \end{equation*}

In the following, we will study numerically the time-harmonic eddy current problems in axisymmetric geometry which can represent correctly the features of our three dimensional problem. 

\section{Axisymmetric domains} \label{section4}
In this section, we choose to consider similar framework and notations introduced in \cite[section 4]{caloz2011influence} in the case of axisymmetric domains and axisymmetric orthoradial data. We suppose that $\Omega_{0}$ and $\Omega_{-}$ are axisymmetric domains with the same axis of rotation denoted by $\rm E_{0}$ which coincides with the \textit{z}-axis. In this case, there exists bi-dimensional "meridian" domains $\Omega^{m},$ $\Omega_{0}^{m},$ and $\Omega_{-}^{m}$ satisfying , in cylindrical coordinates ($r, \theta, z$), the following assumptions 
\begin{equation}
    \begin{array}{l}
    \Omega=\{\mathrm{x} \in \mathbb{R}^{3} / (r, z) \in \Omega^{m}, \theta \in \mathbb{T} \}, \\
    \Omega_{0} = \{\mathrm{x} \in \mathbb{R}^{3} / (r, z) \in \Omega^{m}_{0}, \theta \in \mathbb{T} \}, \\
    \Omega_{-} = \{\mathrm{x} \in \mathbb{R}^{3} / (r, z) \in \Omega^{m}_{-}, \theta \in \mathbb{T}\}.
    \end{array}
\end{equation}
\begin{figure}[h]
    \centering
    \includegraphics[width=0.5\linewidth]{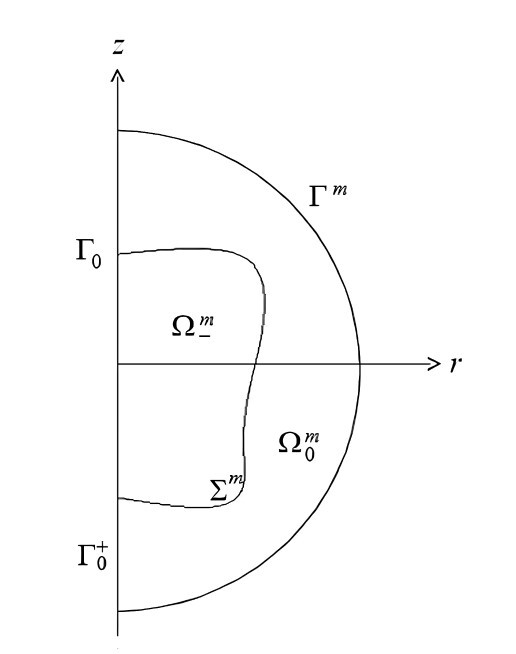}
    \caption{The meridian domain $\Omega=\Omega_{-}^{m} \cup \Omega_{0}^{m} \cup \Sigma^{m}$ with boundary $\partial \Omega^{m}=\Gamma^{m} \cup \Gamma_{0}$}
    \label{Figure 2}
\end{figure}
Here $\mathbb{T} = \mathbb{R}/(2 \pi \mathbb{Z})$ is the one dimensional torus. In Figure \ref{Figure 2},  $\Gamma^{m}$ and $\Sigma^{m}$ are the meridian curves corresponding to $\Gamma$ and $\Sigma,$ and $\Gamma_{0}, \Gamma_{0}^{+}$ are the following subsets of the rotation axis $\mathrm{E}_{0}$
\begin{equation*}
   \begin{array}{lll}
    \Gamma_{0}=\mathrm{E}_{0} \cap \overline{\Omega^{m}} & and & \Gamma_{0}^{+}=\mathrm{E}_{0} \cap \overline{\Omega_{0}^{m}}.
    \end{array}
\end{equation*}

\subsection{Formulation in cylindrical coordinates}
In this part, we recall the cylindrical coordinates of a vector field $\mathcal{A}$ and the $\mathrm{curl}$ operator:
\begin{itemize}
    \item[$\circ$] For a vector field $\mathcal{A}=(\mathcal{A}_{1}, \mathcal{A}_{2}, \mathcal{A}_{3})$ we denote by ($A_{r}, \ A_{\theta}, \ A_{z}$) its cylindrical components such that 
\begin{equation*}
    \left\{ \begin{array}{l}
        A_{r}(r, \theta, z)= \mathcal{A}_{1}(\mathrm{x}) \cos{\theta} + \mathcal{A}_{2}(\mathrm{x}) \sin{\theta},   \\
        A_{\theta}(r, \theta, z)= -\mathcal{A}_{1}(\mathrm{x}) \sin{\theta} + \mathcal{A}_{2}(\mathrm{x}) \cos{\theta}, \\
        A_{z}(r, \theta, z) = \mathcal{A}_{3}(\mathrm{x}),
    \end{array} \right.
\end{equation*}
and we set $\hat{\mathcal{A}}(r, \theta, z)=(A_{r}(r, \theta, z), A_{\theta}(r, \theta, z), A_{z}(r, \theta, z)).$ 
      \item[$\circ$] The cylindrical components of the $\mathrm{curl}$ operator applied to a vector field $\mathcal{A}$ writes
\begin{equation} \label{curl operator in cylindrical coordinates}
    \left\{ \begin{array}{l}
         (\mathrm{curl} \ \mathcal{A})_{r}= \frac{1}{r} \partial_{\theta} A_{z}-\partial_{z} A_{\theta},\\
         (\mathrm{curl} \ \mathcal{A})_{\theta}=\partial_{z} A_{r}- \partial_{r} A_{z}, \\
         (\mathrm{curl} \ \mathcal{A})_{z}=\frac{1}{z} \partial_{r}(r A_{\theta}),
    \end{array} \right.
\end{equation}
and the divergence operator $\mathrm{div}$ writes 
\begin{equation} \label{divergence in cylindrical coordinates}
    \mathrm{div} \mathcal{A}= \partial_{r} A_{r}+\frac{1}{r} A_{r}+\frac{1}{r} \partial_{\theta} A_{\theta}+ \partial_{z} A_{z}.
\end{equation}
\end{itemize}

\subsection{Axisymmetric orthoradial problem}
\subsubsection{Preliminaries}
For a vector field $\mathcal{A}=(\mathcal{A}_{1}, \mathcal{A}_{2}, \mathcal{A}_{3}),$ we say that 
\begin{itemize}
    \item[$\circ$] $\mathcal{A}$ is \textbf{axisymmetric} if $\hat{\mathcal{A}}$ does not depend on the angular variable $\theta$.
    \item[$\circ$] $\mathcal{A}$ is \textbf{orthoradial} if its components $A_{r}$ and $A_{z}$ are equal to zero.
\end{itemize}

On our axisymmetric configuration, we consider a modification of problem (\ref{mag pot A}) \cite{caloz2011influence}: We take $J_{s}=0$ and impose instead a non-homogeneous boundary condition
\begin{equation}
    \mathcal{A} \times n = G \times \textit{n} \ \ \ \mathrm{on} \ \Gamma,
\end{equation}
for a given smooth data G.
\begin{assumption} \label{Assumption G}
    We assume that G is axisymmetric and orthoradial i.e. 
    \begin{equation} \label{G&g_theta}
    \hat{\mathrm{G}}(r, \theta, z)=(0, g_{\theta}(r, z), 0).
    \end{equation}
\end{assumption}
Under Assumption \ref{Assumption G}, it results that $\mathcal{A}$ is also axisymmetric and orthoradial
\begin{equation} \label{A orthoradial-axisymmetric}
    \hat{\mathcal{A}}(r, \theta, z)=(0, A_{\theta}(r, z), 0),
\end{equation}
see for instance \cite{caloz2011influence,peron2009modelisation} for the proof of similar works.

In that follows, we will drop the notation $\theta$ in $A_{\theta},$ and we will concentrate our asymptotic analysis on this orthoradial component. For the sake of clarity, we consider the following notations.
\begin{notation} \label{notation 4}
We denote by $(n_{1},\ n_{2}, \ n_{3})$ the cartesian coordinates of the unit normal vector \textit{n} on $\Sigma$ inwardly oriented to $\Omega_{-}.$ Since $\Omega_{-}$ is an axisymmetric domain, then it results that $n_{\theta}=0$ and the unit normal vector in cylindrical coordinates writes $\hat{n}=(n_{r}, 0, n_{z})$ \cite[page 166]{peron2009modelisation}. Considering an axisymmetric and orthoradial solution (\ref{A orthoradial-axisymmetric}), we introduce then the orthoradial component of the $\rm curl curl$ operator and the boundary operator $\mathrm{curl} \times \text{n}$ respectively as follows:
\begin{equation} \label{operators in cylindrical coordinates}
    \begin{array}{l}
    \mathrm{D}(r, z; \partial_{r}, \partial_{z})=\partial^{2}_{r}+\frac{1}{r} \partial_{r}+\partial^{2}_{z}-\frac{1}{r^{2}},      \\
    \mathrm{B}(r, z; \partial_{r}, \partial_{z})=n_{r}(\partial_{r}+\frac{1}{r})+n_{z} \partial_{z}\textcolor{blue}{,}      
    \end{array}
\end{equation}
and the divergence operator is free in this case, see for instance (\ref{divergence in cylindrical coordinates}).
\end{notation}

\subsubsection{Variational problem}
By using the change of variables from cartesian to cylindrical coordinates, we associate the following weighted Sobolev space in order to define the orthoradial component \textit{A}(\textit{r},\textit{z}) \cite{bernardi1999spectral,caloz2011influence} 
\begin{equation*}
    \mathrm{V}^{1}_{1,\Gamma^{m}}(\Omega^{m})=\{ \mathrm{u} \in \mathrm{H}^{1}_{1}(\Omega^{m}) \ | \ \mathrm{u} \in \mathrm{L}^{2}_{-1}(\Omega^{m}) \ \ \mathrm{and} \ \ \mathrm{u}=0 \ \ \mathrm{on} \ \ \Gamma^{m}\}.
\end{equation*}
Here, 
\begin{equation*}
    \mathrm{H}^{1}_{1}(\Omega^{m})=\{ \mathrm{u} \in \mathrm{L}^{2}_{1}(\Omega^{m}) \  | \  \partial^{j}_{r} \partial_{z}^{1-j} \mathrm{u} \in \mathrm{L}^{2}_{1}(\Omega^{m}),\ j=0,1\},
\end{equation*}
and for all $\alpha \in \mathbb{R},$ the space $\mathrm{L}^{2}_{\alpha}(\Omega^{m})$ is the set of measurable functions $\mathrm{u}(r, z)$ such that
\begin{equation} \label{norm L21}
    \lVert \mathrm{u} \rVert^{2}_{\mathrm{L}^{2}_{\alpha}(\Omega^{m})}=\int_{\Omega^{m}} |\mathrm{u}|^{2} \ r^{\alpha} dr dz < +\infty.
\end{equation}
The following remark incorporates an essential boundary condition.
\begin{remark} \label{remark 1}
    All functions on $\mathrm{V}^{1}_{1,\Gamma^{m}}$ have null trace on $\Gamma_{0},$ see \cite[Remark $\mathrm{II}$.1.1]{bernardi1999spectral}, and \cite[Remark 4.1]{caloz2011influence} for more details about the proof.
\end{remark}
As a result, we solve the following two-dimensional scalar problem set in $\Omega^{m}.$\\
\textit{Find} $A \in \mathrm{V}^{1}_{1,\Gamma^{m}}(\Omega^{m})+g$ \textit{such that for all} $\mathrm{v} \in \mathrm{V}^{1}_{1,\Gamma^{m}}(\Omega^{m}),$
\begin{equation}
    a(A,\mathrm{v})=0,
\end{equation}
where
\begin{equation*}
   \begin{array}{lll}
     a(\mathrm{u},\mathrm{v}) & = & \displaystyle\nu^{2} \int_{\Omega_{-}^{m}} \Big( \partial_{z} \mathrm{u}^{-} \partial_{z} \overline{\mathrm{v}^{-}}+\dfrac{1}{r}\partial_{r}(r\mathrm{u}^{-})\dfrac{1}{r}\partial_{r}(r\overline{\mathrm{v}^{-}}) \Big) \ r dr dz -\mathrm{i} \omega \sigma \mu_{0} \displaystyle \displaystyle\int_{\Omega_{-}^{m}} \mathrm{u}^{-} \overline{\mathrm{v}^{-}} \  r dr dz\\ & & + \int_{\Omega_{0}^{m}} \Big( \partial_{z} \mathrm{u}^{+} \partial_{z} \overline{\mathrm{v}^{+}}+\dfrac{1}{r}\partial_{r}(r\mathrm{u}^{+})\dfrac{1}{r}\partial_{r}(r\overline{\mathrm{v}^{+}}) \Big) \ r dr dz,
     \end{array}
\end{equation*}
and recalling that $\nu=\dfrac{1}{\sqrt{\mu_{r}}}.$

\subsubsection{Strong form of equations}
According to (\ref{operators in cylindrical coordinates}) and Remark \ref{remark 1}, the orthoradial component $A=(A^{+},A^{-})$ satisfies the following problem
\begin{equation} \label{problem Atheta general}
    \left\{ \begin{array}{llll}
    \mathrm{D}A^{+}=0 & & \mathrm{in} & \Omega_{0}^{m} \\
    \mathrm{D}A^{-}-2 \mathrm{i} \delta^{-2} A^{-}=0  & & \mathrm{in} & \Omega_{-}^{m} \\
    \mathrm{B}A^{+}=\varepsilon \delta \mathrm{B}A^{-} & & \mathrm{on} & \Sigma^{m}, \\
    A^{+}=A^{-} & & \mathrm{on} & \Sigma^{m}, \\
    A^{+}=g_{\theta} & & \mathrm{on} & \Gamma^{m}\cup \Gamma_{0}^{+},
    \end{array} \right.
\end{equation}
where $g_{\theta}$ is defined in Eq. (\ref{G&g_theta}), see for instance \cite[Chapter 8]{peron2009modelisation} for similar work. Under the assumption of orthoradial and axisymmetric data and from notation \ref{notation 4}, we deduce directly that the gauge conditions (\ref{gauge conditions}) in the cylindrical coordinates are satisfied.

\subsection{Multiscale expansion} \label{subsection Multiscale expansion}
In this part, we aim to expand the orthoradial component \textit{A} of the magnetic vector potential using a multiscale expansion. First, we introduce the following geometrical notations. 

\begin{notation}[Geometrical setting]
We set $\xi \mapsto \tau(\xi)=(r(\xi), z(\xi))$ a $\mathcal{C}^{\infty}$ function, $\xi \in (0,  L)$ be an \textit{arc-length coordinate} on the interface $\Sigma^{m},$ and \textit{L} is the length of the curve $\Sigma^{m}.$ Let ($\xi$, \textit{h}) be the associate normal coordinate system in a tubular neighborhood $\mathscr{U}_{-}^{m}$ of $\Sigma^{m}$ inside $\Omega_{-}^{m}$ (see for instance Figure \ref{Tubular neighborhood of interface curve}). Then the normal vector $n(\xi)$ at the point $\tau(\xi)$ can be written as (Frenet frame)
\begin{equation} \label{unit normal vector}
    \begin{array}{llll}
    n(\xi)=(-z'(\xi), r'(\xi)), & &\\
    \end{array}
\end{equation}
where $z'(\xi)=\dfrac{dz}{d\xi},$ and $r'(\xi)=\dfrac{dr}{d\xi}.$
Further, we denote by $k(\xi)$ the curvature of $\Sigma^{m}$ at $\tau(\xi)$ which is defined as \cite{caloz2011influence}
$$ k(\xi)=(r' z'' - z' r'')(\xi).$$  
Finally, we set $\chi$ a smooth cut-off function with support in $\overline{\mathscr{U}_{-}^{m}},$ and equals to 1 in a smaller neighborhood of $\Sigma^{m}$.
\end{notation}
Now, we exhibit expansions series for $A$ which we denote by $A^{+}$ in the dielectric part $\Omega_{0}^{m},$ and by $A^{-}$ in the conducting part $\Omega_{-}^{m}:$

\begin{equation} \label{Asymptotic expansion Atheta}
    \begin{array}{l}
        A^{+}(r, z)= A_{0}^{+}(r, z) + \dfrac{\varepsilon}{\hat{\alpha}} A_{1}^{+}(r, z)+ \mathcal{O}(\varepsilon^{2}),    \\
        A^{-}(r,z)=A^{-}_{0}(r,z;\delta)+\delta A^{-}_{1}(r,z;\delta)+\mathcal{O}(\delta^{2}) \\
\ \ \ \ \ \ \ \ with \ A^{-}_{j}(r,z;\delta)=\chi(h) \  \mathfrak{A}_{j}(\xi,\dfrac{h}{\delta}).
    \end{array}
\end{equation}
Here $\hat{\alpha}=\dfrac{1-i}{2}$ where \textit{i} is the unit complex number. Besides, the profiles $\mathfrak{A}_{j}$ are defined on $\Sigma^{m} \times (0, +\infty),$ and $\mathfrak{A}_{j} \longrightarrow 0$ as $Y_{3}=\frac{h}{\delta} \longrightarrow +\infty$. The symbol $\mathcal{O}(\varepsilon^{2})$ (resp. $\mathcal{O}(\delta^{2})$) means that the remainder is uniformly bounded by $\varepsilon^{2}$ (resp. $\delta^{2}$). Hereafter, we focus on the first terms $A_{0}^{+},$ $\mathfrak{A}_{0},$ $A_{1}^{+}$ and $\mathfrak{A}_{1}.$ \\

\subsubsection{First terms of the asymptotic expansion:} 
We construct the first asymptotics $(A_{0}^{+}, \mathfrak{A}_{0})$ and $(A_{1}^{+}, \mathfrak{A}_{1})$ recursively. Elements of formal derivations are given in appendix \ref{Appendix A}.

First, $A_{0}^{+}$ solves the following problem
\begin{equation} \label{A0+ prob}
    \left\{
    \begin{array}{llll}
    \mathrm{D}A_{0}^{+}=0  & & \mathrm{in} & \Omega_{0}^{m}, \\
    \mathrm{B}A_{0}^{+}=0 & & \mathrm{on} & \Sigma^{m}, \\
    A_{0}^{+}=g_{\theta} & & \mathrm{on} & \Gamma^{m} \cup \Gamma_{0}^{+}.
    \end{array} \right.
\end{equation}
Then the first profile $\mathfrak{A}_{0}$ is defined as follows:
\begin{equation}
    \mathfrak{A}_{0}(\xi, Y_{3}) = A_{0}^{+}(\tau(\xi)) \ \mathrm{e}^{-\frac{Y_{3}}{\hat{\alpha}}},
\end{equation}
where $(\xi,\ Y_{3}) \in \mathbb{T}_{L} \times (0,+\infty),$ noting that $\mathbb{T}_{L}=\mathbb{R}/L \mathbb{Z}.$

The next asymptotic solves the problem below:
\begin{equation} \label{A1+ problem}
    \left\{
    \begin{array}{llll}
    \mathrm{D}A_{1}^{+}=0  & & \mathrm{in} & \Omega_{0}^{m}, \\
    \mathrm{B}A_{1}^{+}=-A_{ 0}^{+} & & \mathrm{on} & \Sigma^{m}, \\
    A_{ 1}^{+}=0 & & \mathrm{on} & \Gamma^{m} \cup \Gamma_{0}^{+}.
    \end{array} \right.
\end{equation}
The second profile $\mathfrak{A}_{1}$ satisfies the following equality
\begin{equation}
  \mathfrak{A}_{1}(\xi, Y_{3}) = \bigg[ \frac{1}{\hat{\alpha} \delta_{0}^{2}} A_{1}^{+}(\tau(\xi))-\frac{Y_{3}}{2} (k+\frac{z'}{r})(\xi) A_{0}^{+}(\tau(\xi)) \bigg] \  \mathrm{e}^{-\frac{Y_{3}}{\hat{\alpha}}},
\end{equation}
where $\delta_{0}=\sqrt{\dfrac{2}{\omega \sigma \mu_{0}}}.$ 
Similarly to \cite[Remark 4.2]{caloz2011influence}, we deduce the following remark.
\begin{remark}
Subsequently, we assume that $g_{\theta},$ defined in Eq. (\ref{G&g_theta}), is a real valued function. Thus, the right hand side  of the boundary value problem (\ref{A0+ prob}) is real. Hence, $A_{0}^{+}$ is a real valued function. Similarly, we deduce that $A_{1}^{+}$ is also a real valued function. 
\end{remark}

\begin{remark}
    The models (\ref{A0+ prob}) and (\ref{A1+ problem}) are independent of any physical parameter introduced in this work, which allow us to approach the solution \textit{A} of the problem (\ref{mag pot A}) in the dielectric part with minimal time and memory allocation as well. We provided a deep numerical study about the computational costs in the special issue \cite{abounumerical} for the bi-dimensional case. We will tackle these issues in a forthcoming work for three-dimensional and axisymmetric geometries.
\end{remark}

\subsubsection{Impedance model}
As a by product of the asymptotic expansion, we get a simpler problem then (\ref{mag pot A}) as follows 
\begin{equation} \label{impedance model} 
    \left\{
    \begin{array}{llll}
    \mathrm{D}A_{1}^{\varepsilon}=0  & & \mathrm{in} & \Omega_{0}^{m}, \\
    \mathrm{B}A^{\varepsilon}_{1}+\frac{\varepsilon}{\hat{\alpha}} A^{\varepsilon}_{1}=0  & & \mathrm{on} & \Sigma^{m}, \\
    A^{\varepsilon}_{1}=g_{\theta} & & \mathrm{on} & \Gamma^{m} \cup \Gamma_{0}^{+},
    \end{array} \right.
\end{equation}
where the second condition in (\ref{impedance model}) is the classical Leontovitch condition, see for instance  \cite{leontovich1948approximate}, \cite[section 6.4]{peron2021magnetic}, and \cite[section 3.1]{abounumerical}.
It is well known that the impedance solution has a high accuracy with respect to the solution of the eddy current problems. Accordingly, we will compare numerically our first asymptotic solutions with the latter impedance solution in the next section by using an analytical procedure.

\section{Radial solutions in cylindrical geometry} \label{section 5}
In that follows, our goal is to provide an analytical study of the considered asymptotic models up to the order two in the dielectric part $\Omega_{0}^{m}$. First, we will introduce geometrical and physical assumptions. Then we will give the expressions of the analytical solutions for the global, asymptotic and impedance problems (\ref{problem Atheta general}), (\ref{A0+ prob}), (\ref{A1+ problem}) and (\ref{impedance model}) respectively that are calculated in appendix \ref{Appendix B}. We will assess the accuracy of the resulting asymptotic solutions numerically in the next section.

\subsection{Framework}
We consider a cylindrical geometry: We assume that $\Omega^{m}$ is an infinite cylinder in length consisting of a ferromagnetic material surrounded by a dielectric domain. Let $R_{1}$ be the radius of the interior ferromagnetic cylinder and $R_{2}$ the radius of the domain $\Omega^{m}$. Recall that we denote by $(r, \theta, z)$ the cylindrical coordinate system where the \textit{z}-axis coincides with the axis of the cylinders $\Omega_{-}^{m}$ and $\Omega^{m},$ and by $(\Vec{e}_{r},\  \Vec{e}_{\theta},\  \Vec{e}_{z})$ the local unit vectors in the cylindrical coordinate system. We will assume that the electric current flows in the dielectric domain $\Omega_{0}$ in the $\Vec{e}_{\theta}$ direction and that is uniformly distributed in the $\Vec{e}_{z}$ direction. In order to solve our problem, we impose a Dirichlet condition at $r=R_{2}$: 
\begin{equation} \label{Dirichlet condition}
    A^{+}(r)= \dfrac{k}{r},
\end{equation}
\begin{figure}[!h]
\centering
\begin{minipage}{0.5\textwidth}
    \includegraphics[width=0.7\linewidth]{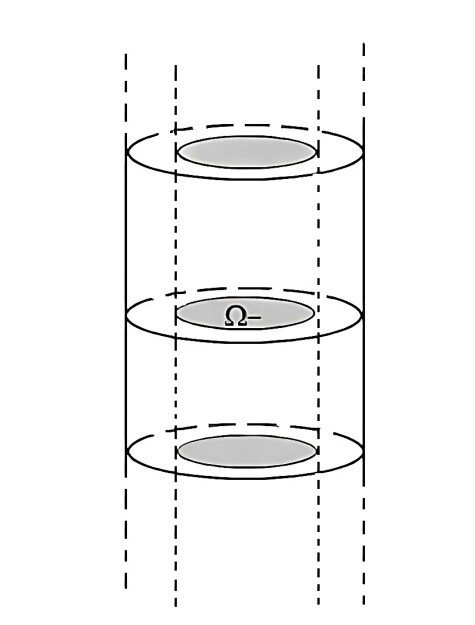}
    \captionsetup{font=small}
    \caption{Considered geometry}
    \label{Geometry}
\end{minipage}
\begin{minipage}[h]{0.49\textwidth}
    \includegraphics[width=0.76\linewidth]{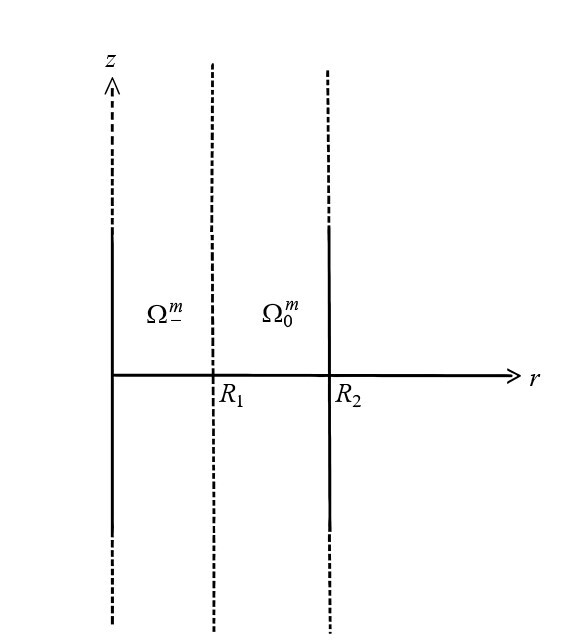}
    \captionsetup{font=small}
    \caption{The meridian domain $\Omega^{m}$ }
    \label{Meridian domains}
    \end{minipage}
\end{figure}
where \textit{k} is a given constant.
The geometry is depicted in Figures \ref{Geometry}-\ref{Meridian domains} below.

\subsection{Analytical solutions}
In the following part, we will exhibit the analytical expressions of $A=(A^{+},A^{-}),$ $A_{0}^{+},$ $A_{1}^{+}$ and $A_{1}^{\varepsilon}$ solutions of problems (\ref{problem Atheta general}), (\ref{A0+ prob}), (\ref{A1+ problem}) and (\ref{impedance model}) respectively.

\subsubsection{Analytical global solutions}
According to \cite[Eq. (A.11)-(A.12)]{bermudez2007transient} and appendix \ref{Appendix B}, the general form of the solutions $A^{+}$ and $A^{-}$ are as follows
\begin{equation} \label{Analytical global solutions}
    \begin{array}{l}
    A^{+}(r)= \dfrac{r}{2} a+ \dfrac{b}{r},       \\
    A^{-}(r)=\mathcal{I}_{1}(\gamma r)  c,    
    \end{array}
\end{equation}
where $\gamma=\sqrt{\omega \sigma \mu_{r} \mu_{0}} e^{+\mathrm{i} \frac{\pi}{4}},$ $\mathcal{I}_{1}$ is the modified Bessel function of the first kind \cite[Chapter 10 - pages 248--250]{olver2010nist}, and \textit{a}, \textit{b} and \textit{c} are constants deduced from the boundary conditions of the above problem (\ref{problem Atheta general}) and having the following expressions:
\begin{equation*}
                \begin{array}{l}
                    a=\dfrac{k}{R_{1}} \dfrac{g_{1}}{g_{2}}, \\
                    b=k-a \dfrac{R_{2}^{2}}{2}, \\
                    c =\big[ \mathcal{I}_{1}^{-1}(\gamma R_{1}) \dfrac{(R_{1}^{2}-R_{2}^{2})}{2R_{1}}\big]a + \dfrac{k}{R_{1}} \mathcal{I}_{1}^{-1}(\gamma R_{1}).
                \end{array}
            \end{equation*}
Noting that $g_{1}$ and $g_{2}$ are constants defined as follows 
\begin{equation*}
    \begin{array}{l}
    g_{1} = \mathcal{I}_{1}(\gamma R_{1})+\gamma R_{1} \mathcal{I}'_{1}(\gamma R_{1}),       \\
    g_{2} = R_{1} \mu_{r} \mathcal{I}_{1}(\gamma R_{1})-g_{1} \dfrac{(R^{2}_{1}-R^{2}_{2})}{2R_{1}}.      
    \end{array}
\end{equation*} 
\subsubsection{Analytical asymptotic solutions} \label{subsection analytical solutions}
We find the analytical expressions of $A_{0}^{+}$ and $A_{1}^{+}$ by a simple integration of the first equations in $\Omega_{0}^{m}$ of (\ref{A0+ prob}) and (\ref{A1+ problem}). The latter asymptotics have the following general form:
\begin{equation} \label{analytical asymptotics}
   \begin{array}{l}
        A_{0}^{+}(r)= \dfrac{r}{2} a_{0} + \dfrac{b_{0}}{r},  \\
        A_{1}^{+}(r)= \dfrac{r}{2} a_{1} + \dfrac{b_{1}}{r},
   \end{array}
\end{equation}
where $(a_{0}, b_{0})$ and $(a_{1}, b_{1})$ are constants that are deduced from the boundary conditions of (\ref{A0+ prob}) and (\ref{A1+ problem}) respectively. Their expressions are given below:
\begin{equation*}
    \begin{array}{lll}
         a_{0}=0  & \mathrm{and} & b_{0}=k, \\
         a_{1}=\dfrac{k}{R_{1}} & \mathrm{and} & b_{1}= -\dfrac{k R_{2}^{2}}{2 R_{1} }.
        \end{array}
\end{equation*}
Finally, the first asymptotic solutions have the following analytical form:
\begin{itemize}
    \item[-] \textit{order 1:} $A_{0}^{+}(r)=\dfrac{k}{r},$
    \item[-] \textit{order 2:} $A_{0}^{+}(r)+\dfrac{\varepsilon}{\hat{\alpha}} A_{1}^{+}(r)=\dfrac{k}{r}+\dfrac{\varepsilon}{\hat{\alpha}}(\dfrac{k}{R_{1}} \dfrac{r}{2} -\dfrac{k R_{2}^{2}}{2 R_{1}} \ \dfrac{1}{r}).$
\end{itemize}
\subsubsection{Analytical impedance solution}
We find the analytical expressions of $A_{1}^{\varepsilon}$ by a simple integration of the first equations in $\Omega_{0}^{m}$ of the impedance model (\ref{impedance model}). The latter solution has the following general form:
\begin{equation} \label{analytical impedance solution}
   \begin{array}{l}
        A_{1}^{\varepsilon}(r)= \dfrac{r}{2} a_{2}^{\varepsilon}+ \dfrac{b_{2}^{\varepsilon}}{r},  \\
   \end{array}
\end{equation}
where $a_{2}^{\varepsilon}$ and $b_{2}^{\varepsilon}$ are constants that are deduced from the boundary conditions of (\ref{impedance model}). Their expressions are given below:
\begin{equation*}
   \begin{array}{lll}
       a_{2}^{\varepsilon}=\dfrac{\varepsilon}{\hat{\alpha}}\dfrac{2 k}{\zeta^{\varepsilon}} & \mathrm{and} & 
       b_{2}^{\varepsilon}=k-a_{2}^{\varepsilon} \dfrac{R_{2}^{2}}{2}.
                 \end{array} 
             \end{equation*}
where $\zeta^{\varepsilon}=(\dfrac{\varepsilon}{\hat{\alpha}} (R_{2}^{2}-R_{1}^{2})+2R_{1}).$

\section{Numerical results} \label{section6}
This section is devoted to establish numerical experiments concerning the above analytical solutions which have been implemented using Python 3 \cite{hellmann2011python,lutz2013learning}. The considered physical parameters are illustrated in Table \ref{parameters}. We choose arbitrarily small radius $R_{1}$ and $R_{2}$. Moreover, the choice of the constant \textit{k}, introduced in the Dirichlet condition (\ref{Dirichlet condition}), is also arbitrary. 
\begin{table}[!h]
\centering
\begin{tabular}{ ll } 
 \hline
   Parameters & Value \\ 
\hline
   Relative permeability ($\mu_{r}$) & 4000 \\ 
   Conductivity ($\sigma$) & 2E+06 S/m \\
   Frequency (f) & 10 Hz \\
   Inner radius ($ R_{1}$) & 0.03 m \\
   Outer radius ($ R_{2}$) & 0.04 m \\
   Skin depth ($\delta$) & 1.779E-03 m \\
   Epsilon ($ \varepsilon $) & 1.41E-01 \\
   \textit{k} & 1 \\
\hline
\end{tabular}
\caption{Physical and numerical parameters}
\label{parameters}
\end{table}

Our numerical analysis follows the spirit of \cite{abounumerical}, where the finite element method was applied to demonstrate the accuracy of our asymptotic approach in a bi-dimensional setting for the eddy current problems in linear ferromagnetic materials. Indeed, the finite element solution was imposed as the reference solution in \cite{abounumerical}, since its corresponding analytical expression was not obvious to calculate.

First, we study the errors of orders one and two in the interval [$R_{1}$, $R_{2}$] of the dielectric part $\Omega_{0}^{m}$, in order to demonstrate the accuracy of our first asymptotic models.

\begin{figure}[!htb]
     \centering
     \begin{minipage}[h]{0.5\textwidth}
     \centering
         \includegraphics[width=0.96\textwidth]{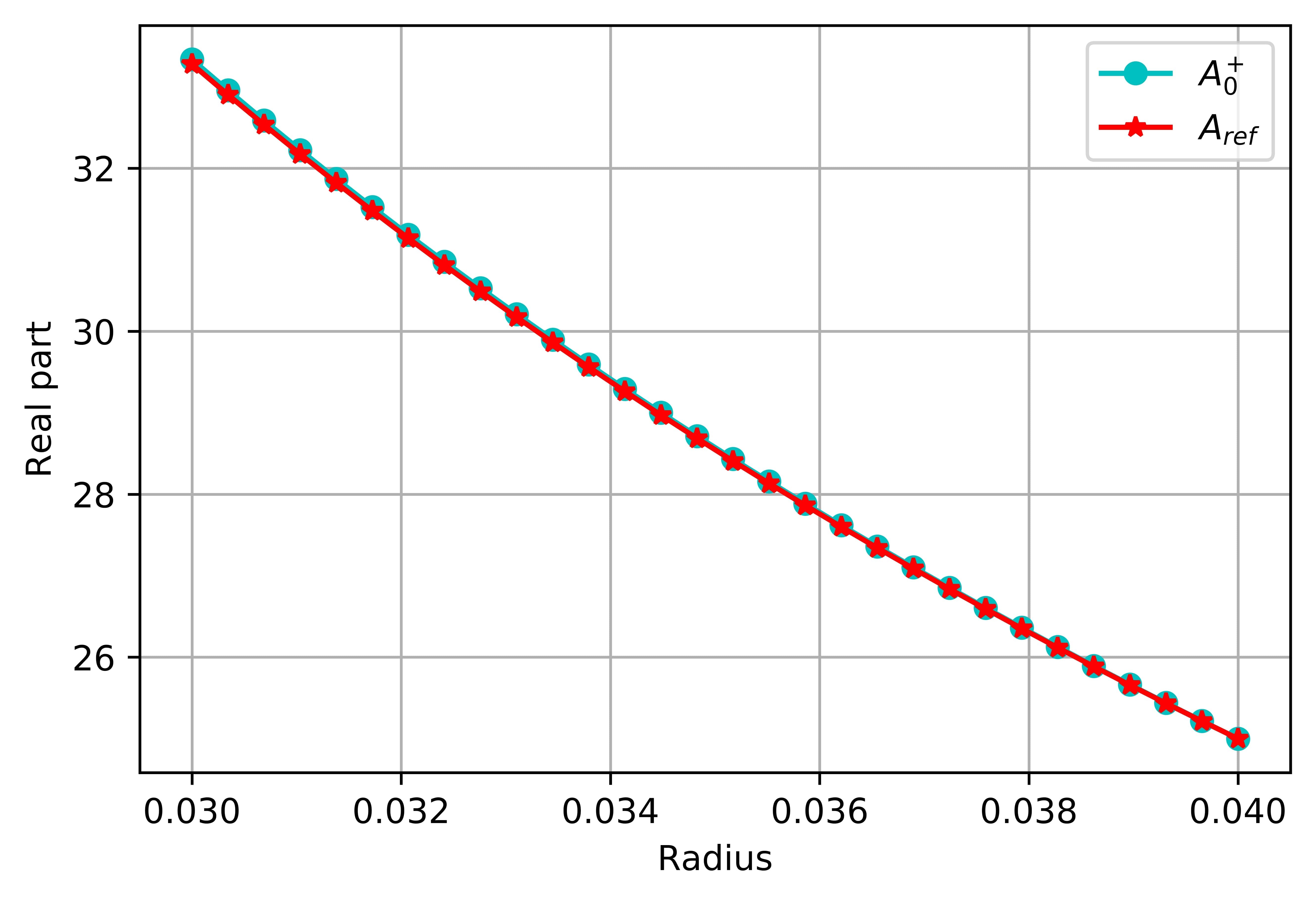}
         \caption{Real parts of $A_{0}^{+}$ and $A_{ref}$}
         \label{Real A0 and A ref}
     \end{minipage} \hfill
     \begin{minipage}[h]{0.49\textwidth}
         \centering
         \includegraphics[width=1\textwidth]{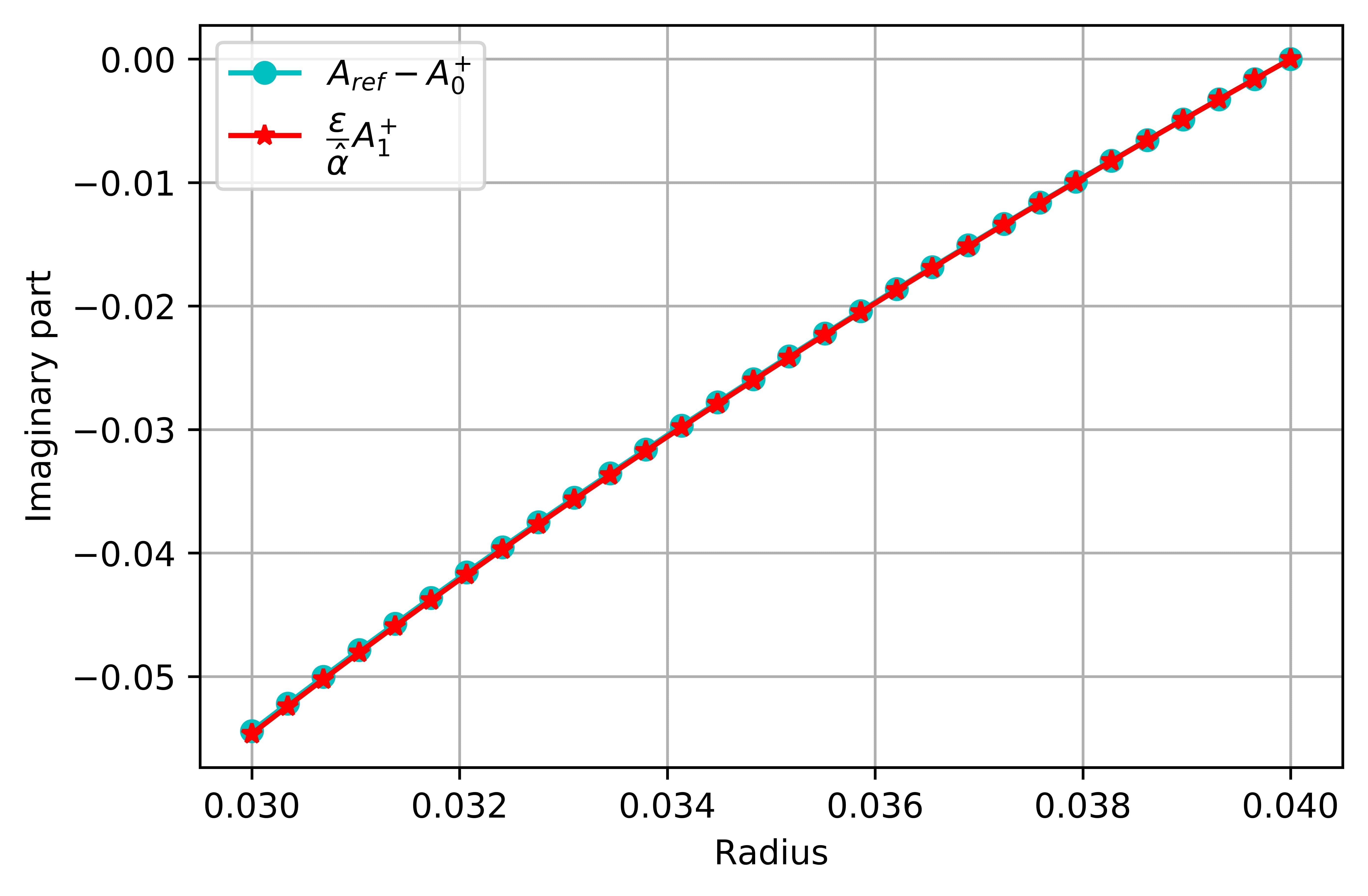}
         \caption{Imaginary parts of $A_{ref}-A_{0}^{+}$ and $\dfrac{\varepsilon}{\hat{\alpha}}A_{1}^{+}$}
         \label{Imag Error 1 and correction 1}
     \end{minipage}
\end{figure}
\begin{figure}[!htb]
     \centering
         \includegraphics[width=0.52\textwidth]{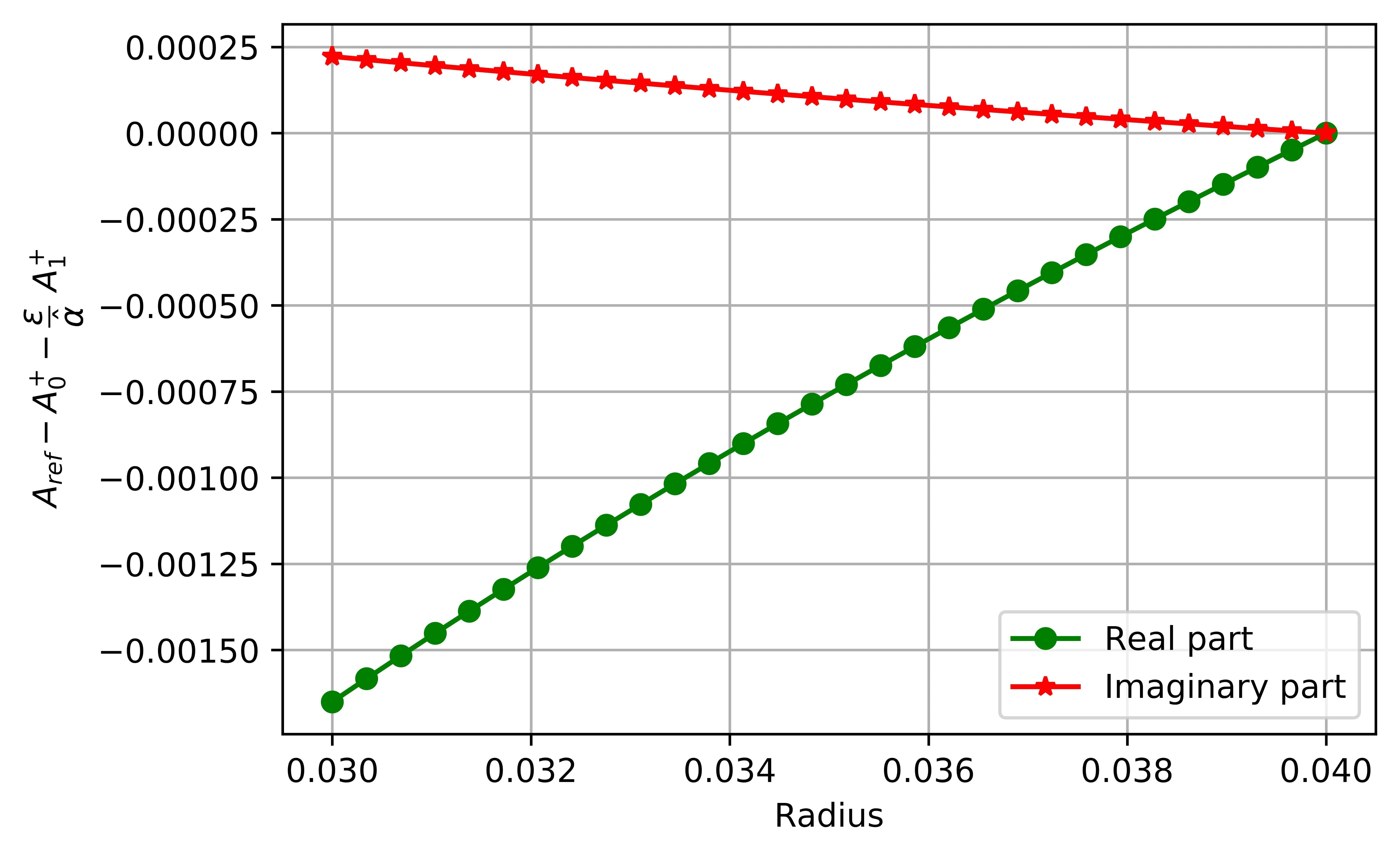}
         \caption{Real and Imaginary parts of the second order error $A_{ref}-A_{0}^{+}-\dfrac{\varepsilon}{\hat{\alpha}}A_{1}^{+}$}
         \label{Real and imaginary parts of the second order model}
\end{figure}

In Figure \ref{Real A0 and A ref}, the real part of the first asymptotic $A_{0}^{+}$ is coherent with that of the reference solution in $\Omega_{0}^{m}$ denoted by $A_{ref}.$ Moreover, we plot the error of the first order $A_{ref}-A_{0}^{+}$ and its correction $\dfrac{\varepsilon}{\hat{\alpha}}A_{1}^{+}$ for the imaginary parts in Figure \ref{Imag Error 1 and correction 1}. We remark that the corresponding graphs are approximately the same. Noting that we have the same results for the real parts. Thus, we conclude that we have a good precision and correction of the first asymptotic solution. 
Since $A_{0}^{+}$ is real unlike $A_{ref}$, then it would be interesting to study the accuracy of the second order solution $A_{0}^{+}+\dfrac{\varepsilon}{\hat{\alpha}} A_{1}^{+}.$ This accuracy is established by the implementation of the second order error $A_{ref}-A_{0}^{+}-\dfrac{\varepsilon}{\hat{\alpha}} A_{1}^{+}.$ Indeed, Fig. \ref{Real and imaginary parts of the second order model} ensures that the real and imaginary parts of the latter error are small enough since they are  less than $1 \% $. As a consequence, we deduce the good approximation of the asymptotic solution $A_{0}^{+} + \dfrac{\varepsilon}{\hat{\alpha}}A_{1}^{+}.$

Next, it is useful in power electronics to study the accuracy of our approach for different values of the physical parameters introduced in our work. In this context, we aim to compare the relative errors of the first asymptotic solutions with that of the impedance solution satisfying problem (\ref{impedance model}). We plot the convergence graphs with the log-log scale defining the following relative $\mathrm{L}^{2}_{1}$ error:
\begin{equation*}
    Error=\dfrac{\Vert A_{ref} - A_{num} \Vert_{\mathrm{L}^{2}_{1}([R_{1},R_{2}])}}{\Vert A_{ref} \Vert_{\mathrm{L}^{2}_{1}([R_{1},R_{2}])}}
\end{equation*}
where the norm $\Vert \cdot \Vert_{\mathrm{L}^{2}_{1}[R_{1},R_{2}]}$ is defined as follows $$\Vert \mathrm{u} \Vert_{\mathrm{L}^{2}_{1}([R_{1},R_{2}])}=\Big(\int_{[R_{1},R_{2}]} |\mathrm{u}|^{2} \ r dr \Big)^{\frac{1}{2}}.$$ Besides, $A_{ref}$ is the reference solution in $\Omega_{0}^{m},$ and $A_{num}$ is the first asymptotic model $A_{0}^{+},$ the second asymptotic model $A_{0}^{+}+\dfrac{\varepsilon}{\hat{\alpha}} A_{1}^{+},$ or the impedance solution $A_{1}^{\varepsilon}.$ We shall consider the case of varying relative magnetic permeabilities and frequencies as well. The geometry and the considered physical parameters are depicted in Figures \ref{Geometry}-\ref{Meridian domains}, and Table \ref{parameters}. 
\begin{figure}[!htb]
     \centering
     \begin{minipage}[h]{0.5\textwidth}
         \includegraphics[width=1\textwidth]{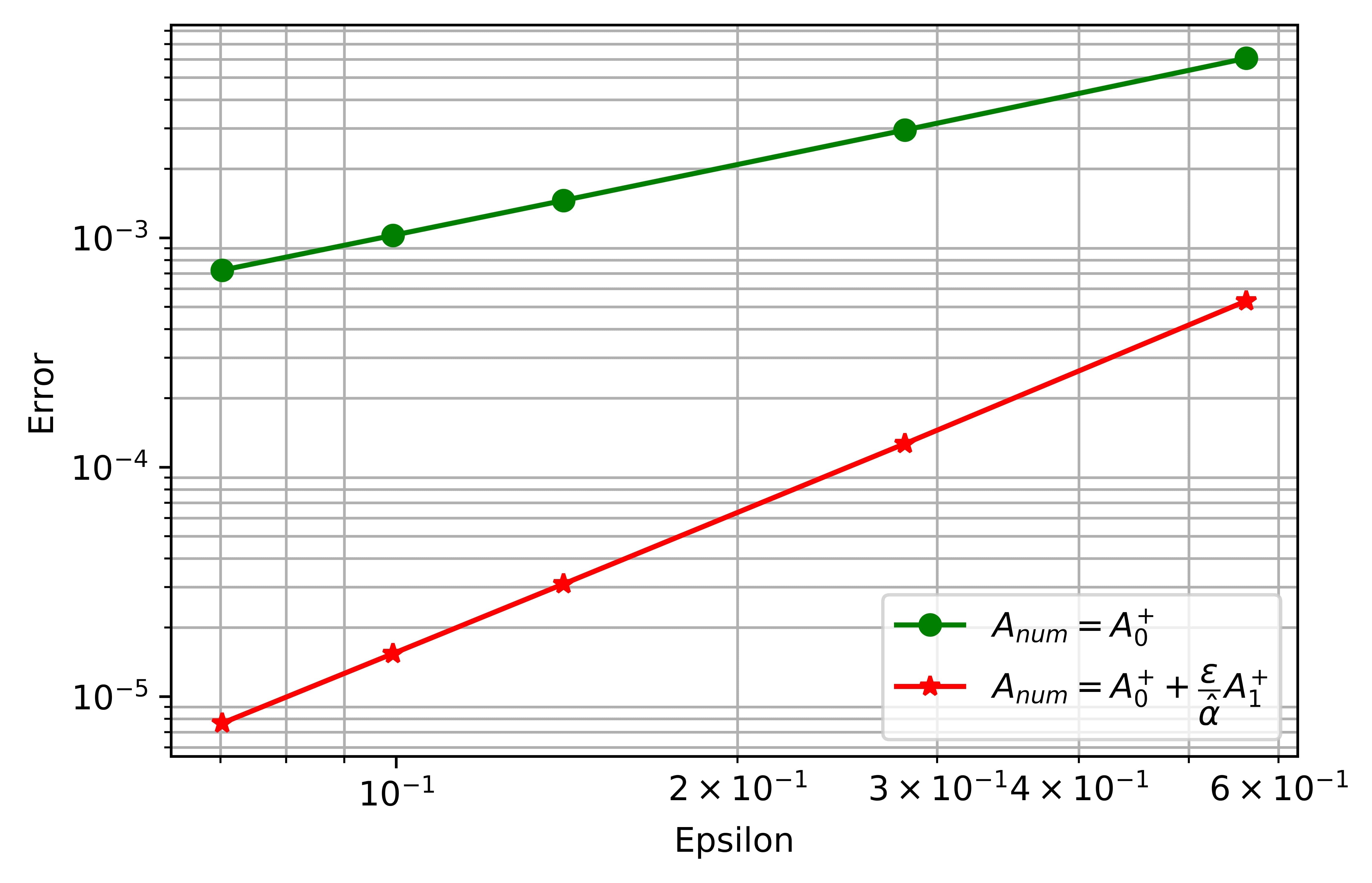}
         \caption{Relative $\mathrm{L}^{2}_{1}$ error of the asymptotic solutions of order 1 and 2 for f = 10 Hz}
         \label{Relative errors order 1 and 2}
     \end{minipage} \hfill
     \begin{minipage}[h]{0.49\textwidth}
         \includegraphics[width=1\textwidth]{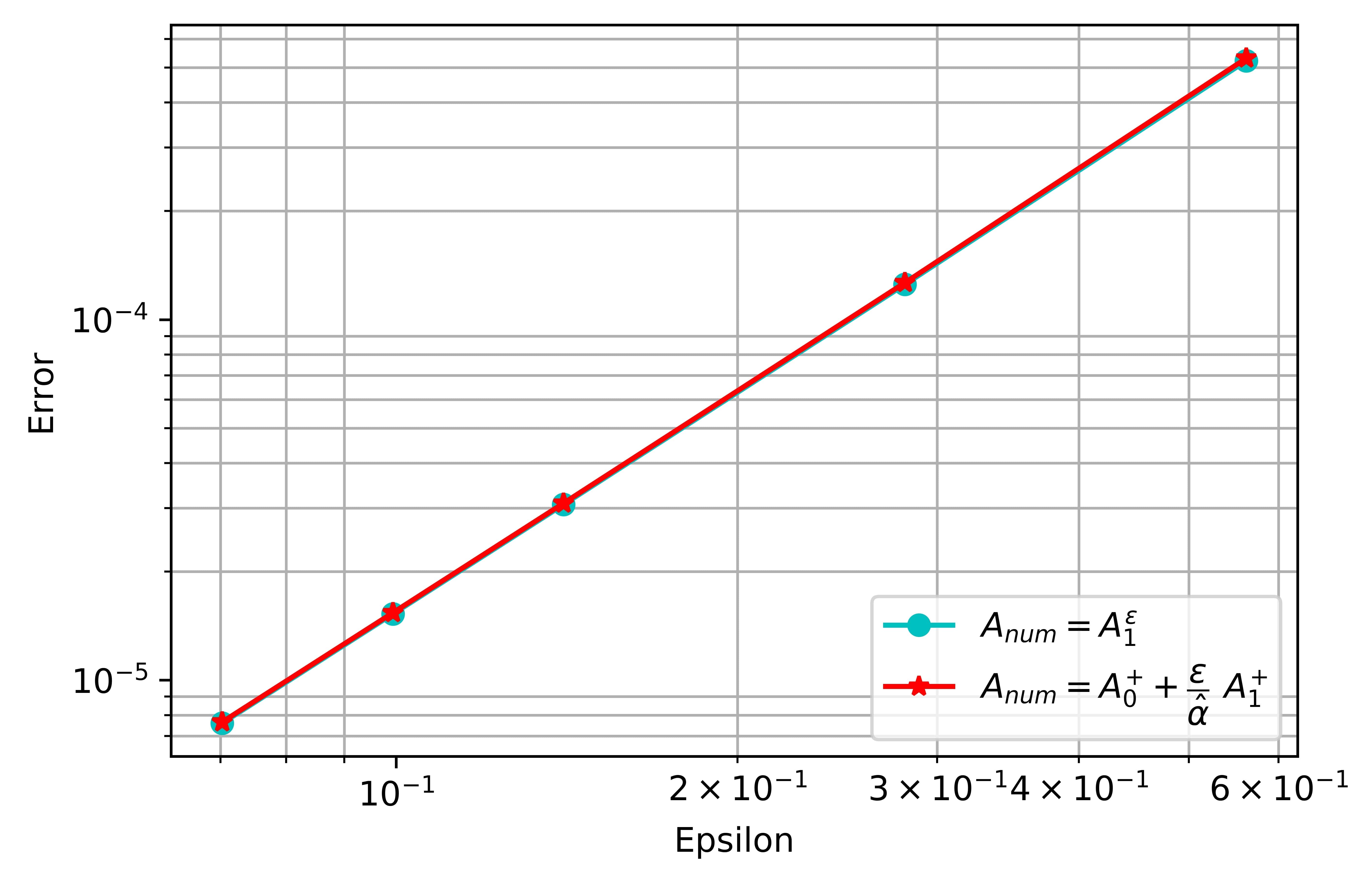}
         \caption{Relative $\mathrm{L}^{2}_{1}$ error of the impedance solution and the second order model for f = 10 Hz}
         \label{Relative errors of the impedance and second order solution}
     \end{minipage}
\end{figure}
We recall that our work is restricted to eddy currents in linear ferromagnetic materials. In this regard, we suppose in Figures \ref{Relative errors order 1 and 2}-\ref{Relative errors of the impedance and second order solution} that the relative magnetic permeabilities are high and between 250 and 16000, and the frequency is 10 \text{Hz}. Noting that this latter range was chosen arbitrarily. We ensure in Figure \ref{Relative errors order 1 and 2} that when the small parameter $\varepsilon$ decreases, the convergence rate of the relative $\mathrm{L}^{2}_{1}$ error is of order 1 for the model $A_{0}^{+}$ and of order 2 for $A_{0}^{+}+\dfrac{\varepsilon}{\hat{\alpha}} A_{1}^{+}$. Moreover, the asymptotic solution $A_{0}^{+}+\dfrac{\varepsilon}{\hat{\alpha}} A_{1}^{+}$ provides an approximation to the reference solution which is of the same rate as the impedance solution $A_{1}^{\varepsilon},$ since the corresponding errors exhibited in Figure \ref{Relative errors of the impedance and second order solution} behave in a similar manner with the variation of $\varepsilon.$
\begin{figure}[!htb]
     \centering
     \begin{minipage}[h]{0.5\textwidth}
         \includegraphics[width=1\textwidth]{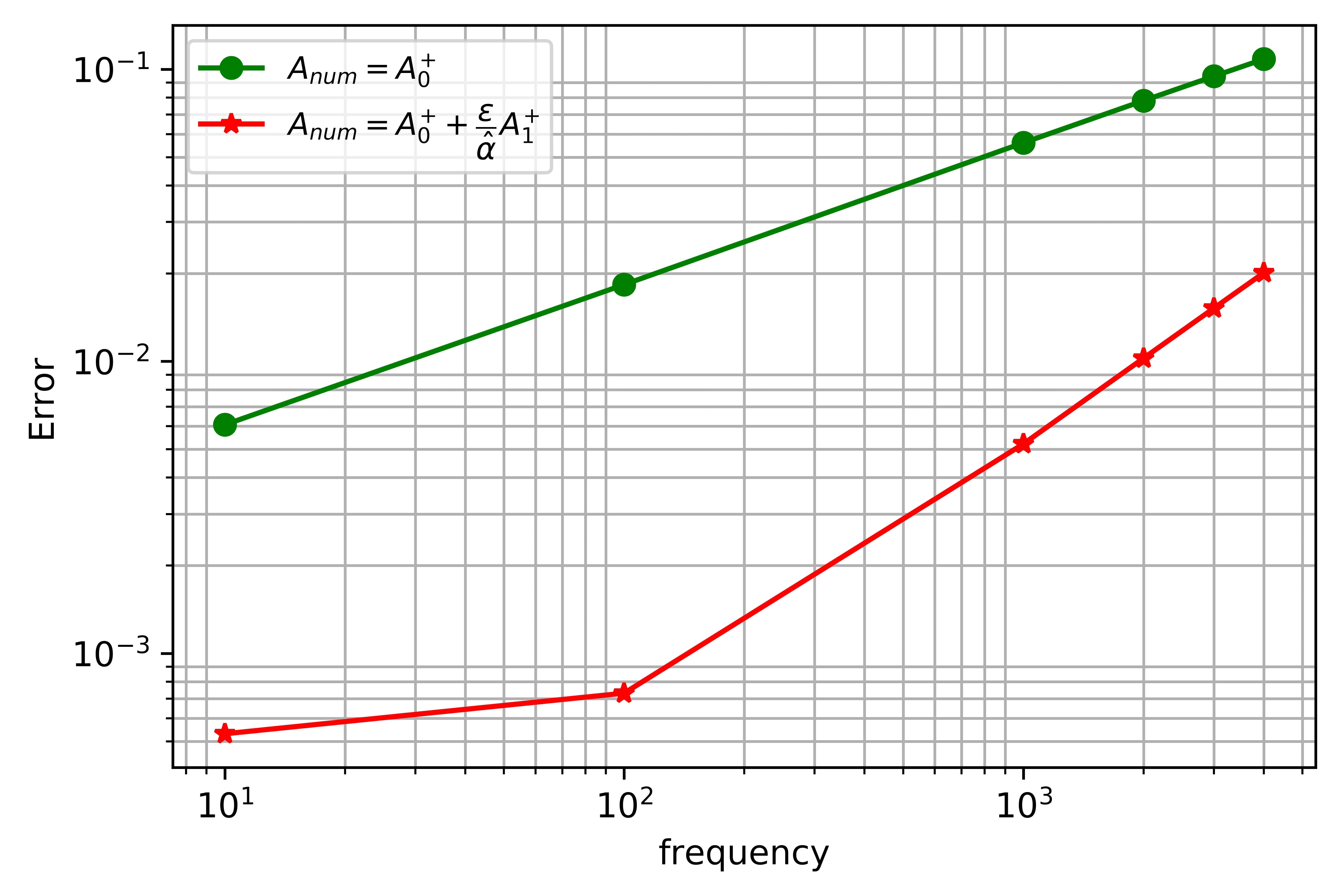}
         \caption{Relative $\mathrm{L}^{2}_{1}$ error versus frequency for the asymptotic models for \\ $\mu_{r} = 250$}
         \label{Relative infinite error versus frequency for mur=250}
     \end{minipage}
     \begin{minipage}[h]{0.49\textwidth}
         \includegraphics[width=1\textwidth]{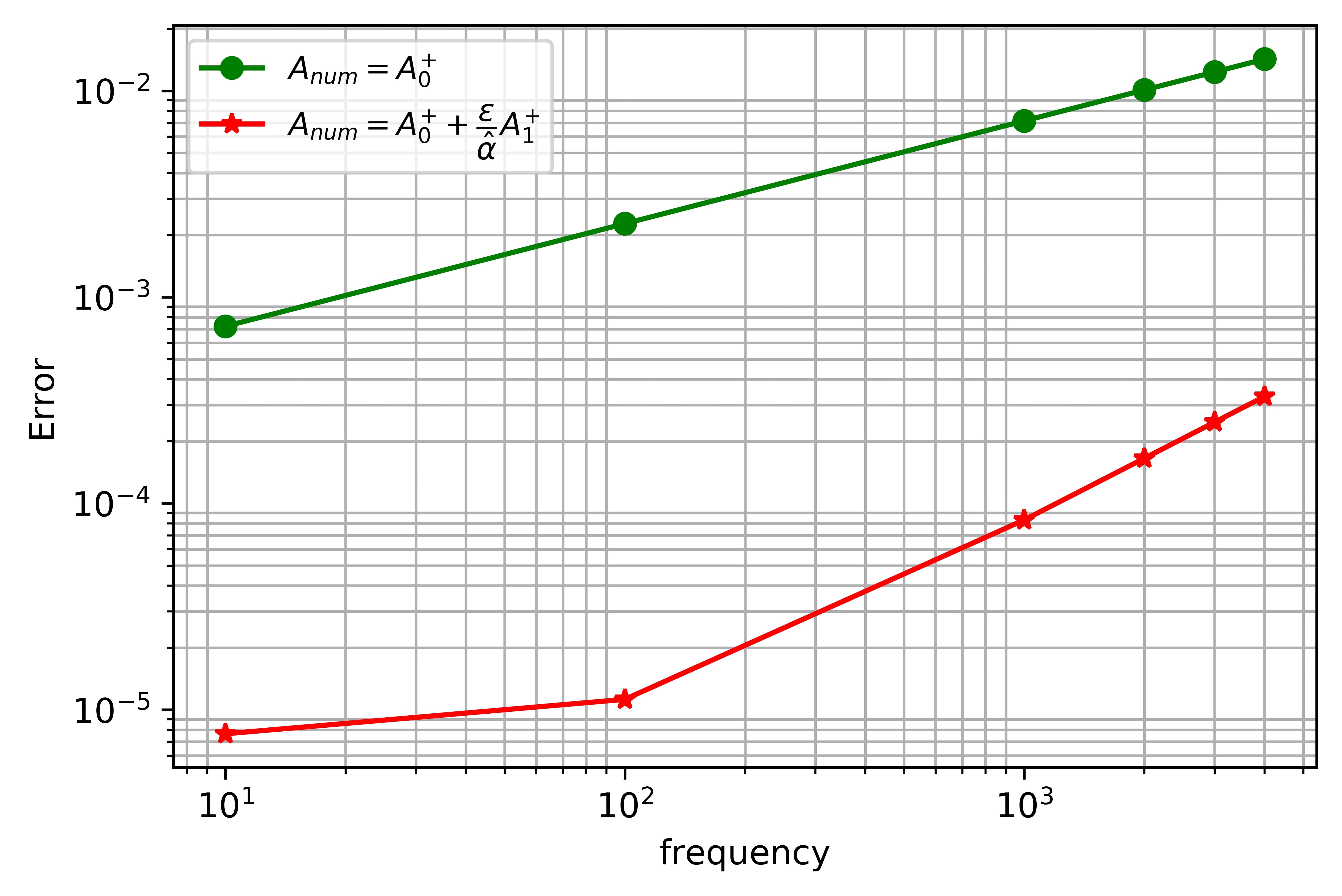}
         \caption{Relative $\mathrm{L}^{2}_{1}$ error versus frequency for the asymptotic models for \\ $\mu_{r} = 16000$}
         \label{Relative infinite error versus frequency for the asymptotic models for mur16000}
     \end{minipage}
\end{figure}

Now, the relative $\mathrm{L}^{2}_{1}$ errors versus frequency for the asymptotic solutions are depicted in Figures \ref{Relative infinite error versus frequency for mur=250}-\ref{Relative infinite error versus frequency for the asymptotic models for mur16000} and for two different relative permeabilities. It is important to recall that our asymptotic approach is accomplished when the parameter $\varepsilon=\dfrac{1}{\mu_{r} \delta}$ is small i.e. less than one. In order to achieve this assumption, we must have the frequency between 10 Hz and 30 Hz, for $\mu_{r}=250,$ and between 10 Hz and 2 kHz for $\mu_{r}=16000.$ We will focus our analysis on these ranges of frequencies. On the one hand, when $\mu_{r}=250$, the first order relative $\mathrm{L}^{2}_{1}$ error is less than 1$\%$ for a range of frequencies between 10 Hz and 25 Hz. The range of frequencies becomes larger between 10 Hz and 1.6 kHz, when the relative permeability increases to 16000. On the other hand, we get slightly better results when the order of the asymptotic model is two. Indeed, the second order relative $\mathrm{L}^{2}_{1}$ error is less than $1\%$ for f in $[10, 30],$ when $\mu_{r}=250,$ and for a wider range of frequencies $[10,2000]$ when $\mu_{r}=16000.$ It is worthwhile to note here that in both cases we remark a slope break for the graph corresponding to the second order error. Precisely, the slope is less important for low frequencies f in $[ 10,100]$, so in the case of not enough "small" skin depth $\delta$ and when $\varepsilon$ tends to zero, than that beyond this latter range of frequencies. As a future work, it could be interesting to study the solution of higher order, for instance from 30 Hz on when $\mu_{r} \geq 250,$ and from 2 kHz on when $\mu_{r} \geq 16000$ in order to get a relative $\mathrm{L}^{2}_{1}$ error less than $1\%$ for a deeper range of frequencies and when $\varepsilon$ must be small as well. 

\section{Conclusion and perspectives} \label{section7}
In conclusion, this paper provides efficient asymptotic models for axisymmetric eddy current problems in linear ferromagnetic materials. Our numerical experiments, established analytically for a special class of unbounded domains, confirm that the proposed asymptotic approach with two components suffices to ensure high accuracy for the case of low frequencies.

As a future work, we aim to study analytically and numerically the case of smooth and bounded geometries. Moreover, the multiscale approach for the three-dimensional eddy current problem as well as proofs of error estimates will be tackled in the PhD thesis \cite{abounum} which is in preparation. It would be useful to expand our analysis numerically for the 3D case. Finally, we recall that our approach does not fit near edges and corners on the conductor interface. In this perspective, we will investigate in a forthcoming work an asymptotic procedure that provides reduced computational costs concerning geometrical singularities of the eddy current problems in linear ferromagnetic materials.

\appendix
\section{Elements of derivation for the multiscale expansion} \label{Appendix A}
In this section, we will derive the terms of the asymptotic expansions introduced in (\ref{Asymptotic expansion Atheta}) at any order $n \in \mathbb{N}$ as well as their governing equations having in mind that the orthoradial component of the magnetic vector potential $A=(A^{+}, A^{-})$ satisfies the following problem 
 
\begin{equation} \label{Appendix problem Atheta general}
    \left\{ \begin{array}{llll}
    \mathrm{D}A^{+}=0 & & \mathrm{in} & \Omega_{0}^{m} \\
    \mathrm{D}A^{-}-2 \mathrm{i} \delta^{-2} A^{-}=0  & & \mathrm{in} & \Omega_{-}^{m} \\
    \mathrm{B}A^{+}=\varepsilon \delta \mathrm{B}A^{-} & & \mathrm{on} & \Sigma^{m}, \\
    A^{+}=A^{-} & & \mathrm{on} & \Sigma^{m}, \\
    A^{+}=g_{\theta} & & \mathrm{on} & \Gamma^{m} \cup \Gamma_{0}^{+},
    \end{array} \right.
\end{equation}
and recalling that from notation \ref{notation 4}, we deduce directly that the gauge conditions (\ref{gauge conditions}) in the cylindrical coordinates are satisfied.
We remind that the magnetic potential $A$ is concentrated on the boundary $\Sigma^{m}$ and decays rapidly inside the conductor. Hence, it is convenient to use a local "normal coordinate system" in a tubular neighborhood $\mathscr{U}_{-}^{m}$ of $\Sigma^{m}$ inside $\Omega_{-}^{m}.$ 

We denote by ($\Vec{e}_{r}$, $\Vec{e}_{\theta}$, $\Vec{e}_{z}$) the basis associated with the cylindrical coordinates $(r, \theta, z)$. In the basis ($\Vec{e}_{r}$, $\Vec{e}_{z}$), recall that $(r(\xi), \ z(\xi))=\tau(\xi), \ \xi \in (0, L)$ is an arc length coordinate on the interface $\Sigma^{m},$ and ($\xi$, $h$) is the associate normal coordinate system. The normal vector $n(\xi)$ at the point $\tau(\xi)$ writes $$n(\xi)=(-z'(\xi), \ r'(\xi)).$$ 
Hence, the tubular neighborhood $\mathscr{U}_{-}^{m}$ of $\Sigma^{m}$ inside $\Omega^{m}_{-}$ is represented by the parameterization below 
$$ \mathscr{U}^{m}_{-}= \Psi(\mathbb{T}_{L} \times [0, \ h_{0})),$$
where $\Psi$ is the change of coordinates defined by $$\Psi: (\xi, h) \mapsto (r, z),$$ noting that $r=r(\xi)-h z'(\xi),$ and $z=z(\xi)+h r'(\xi)$. 
Recalling that the curvature is defined by the following equality
$$k(\xi)=(r' z'' - z' r'')(\xi).$$
It is important to note that for $h_{0} < \dfrac{1}{\Vert k \Vert}_{\infty},$ the change of coordinates $\Psi$ is
a $\mathcal{C}^{\infty}$-diffeomorphism from the cylinder $\mathbb{T}_{L} \times [0, \ h_{0})$ into $\mathscr{U}_{-}^{m}$. In contrast, when the radius of the interface curvature is very small, for example less then \textit{h}, we get a rounded corner on $\Sigma^{m},$ and at the limit i.e. when the curvature is infinite we get a sharp corner (see for instance Remark 2 in \cite{abounumerical}). In this paper, our interest lies in the case where the interface $\Sigma^{m}$ is smooth. The latter cases are beyond the scope of this research. 

\begin{figure}[t]
    \centering
    \includegraphics[width=0.35\linewidth]{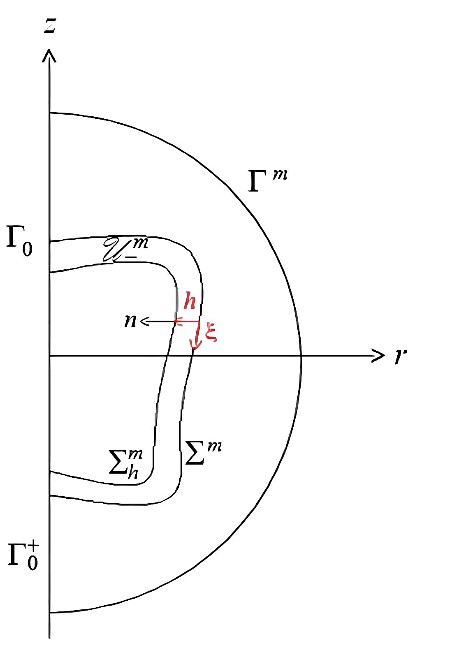}
    \caption{A tubular neighborhood of $\Sigma^{m}$ }
    \label{Tubular neighborhood of interface curve}
\end{figure}

In the following section, we aim to identify the profiles $\mathfrak{A}_{n}$ as well as the asymptotic models $A_{n}^{+}$ introduced in (\ref{Asymptotic expansion Atheta}) at any order $n \in \mathbb{N}.$ To do that, we first expand the interior operator D and the boundary operator B in power series of the skin depth $\delta$ by applying the change of variables $\Psi$ and the scaling $Y_{3}=\dfrac{h}{\delta}.$ Then we plug the resulting expressions in the problem (\ref{Appendix problem Atheta general}). Finally, by identifying with the same power in $\varepsilon,$ $\delta,$ and $\mu_{r},$ we get the coefficients of the asymptotic expansion satisfying the family of boundary value problems at any order $n \in \mathbb{N}$. For the sake of simplicity, we will explicit the first asymptotics $\mathfrak{A}_{n}$ and $A_{n}^{+}$ for $n=0,1$. 

\subsection{Expansion of the operators}
In this part, we expand the operators D and B in the same manner as in \cite{caloz2011influence}. 

Performing the change of the scaling $Y_{3}=\dfrac{h}{\delta}$ and the change of variables $\Psi$, the interior operator D writes in coordinates ($\xi$, \textit{h}) as 
\begin{equation}
    \mathrm{D}=\delta^{-2} \Big[ \partial^{2}_{Y_{3}}+ \delta \mathrm{D}_{1}+\delta^{2} \mathrm{R}_{\delta} \Big]
\end{equation}
where $\mathrm{D}_{1}(\xi,Y_{3}; \partial_{\xi}, \partial_{Y_{3}}) = -(k+\frac{z'}{r})(\xi) \partial_{Y_{3}}$ and $R_{\delta}$ is an operator, which has smooth coefficients in $Y_{3}$ and $\xi,$ bounded in $\delta.$
Hence, we get 
\begin{equation}
    \delta^{2} \mathrm{D} = \sum\limits_{n \geq 0} \delta^{n} C_{n}
\end{equation}
where
\begin{equation*}
 \begin{array}{l}
      C_{0}=\partial_{Y_{3}}^{2},   \\
      C_{1}=-(k+\frac{z'}{r})(\xi) \partial_{Y_{3}}. 
 \end{array}
 \end{equation*}
Similarly, there holds $\mathrm{B}=\delta^{-1} \partial_{Y_{3}}-\frac{z'}{r}(\xi)$ on the interface $\Sigma.$

\subsection{Equations of the coefficients of the magnetic potential}
In this section, we define $v_{\delta}(\xi, Y_{3})=A^{-}(\mathrm{x})$ in $\mathscr{U}_{-}^{m}.$ After the scaling $h \longmapsto Y_{3}=\dfrac{h}{\delta}$ in $\mathscr{U}^{m}_{-},$ the problem (\ref{Appendix problem Atheta general}) writes 
\begin{equation} \label{Atheta+}
    \left\{ \begin{array}{lll}
    \mathrm{D} A^{+}=0     & & \mathrm{in}\ \Omega_{0}^{m},  \\
    \mathrm{B} A^{+}=\varepsilon (\partial_{Y_{3}}-\frac{z'}{r}) v_{\delta} & & \mathrm{on} \ \Sigma^{m}, \\
    A^{+}=0 & & \mathrm{on} \ \Gamma^{m} \cup \Gamma_{0}^{+},
    \end{array} \right.
\end{equation}
and, 
\begin{equation} \label{vdelta}
    \left\{ \begin{array}{lll}
    (\partial^{2}_{Y_{3}}-(\frac{1}{\hat{\alpha}})^{2})v_{\delta}-\sum\limits_{n \geq 1} \delta^{n} C_{n} v_{\delta}=0 & & \mathrm{in} \ \mathbb{T}_{L} \times (0, +\infty),       \\
    v_{\delta}=A^{+} & & \mathrm{on} \ \mathbb{T}_{L} \times \{0\}.      
    \end{array} \right.
\end{equation}
Now we plug the ansatz,
\begin{equation*}
    \begin{array}{lll}
    A^{+} \sim \sum\limits_{n \geq 0} \big( \dfrac{\varepsilon}{\hat{\alpha}}\big)^{n} A_{n}^{+}(\mathrm{x}) & & \mathrm{in} \ \Omega_{0}^{m},  \\
    \end{array}
\end{equation*}
and 
\begin{equation*}
    \begin{array}{lll}
    v_{\delta} \sim \sum\limits_{n \geq 0} \delta^{n} \mathfrak{A}_{n}(\xi, Y_{3}) & & \mathrm{in} \ \mathscr{U}_{-}^{m},  \\
    \end{array}
\end{equation*}
with $\mathfrak{A}_{n} \longrightarrow 0$ as $Y_{3} \longrightarrow + \infty$ in (\ref{Atheta+}) and (\ref{vdelta}). Then by identification of terms in power of $\varepsilon,$ $\delta$ and $\mu_{r},$ the profiles $\mathfrak{A}_{n}$ and $A_{n}^{+}$ satisfy the family of problems coupled by their conditions on the interface $\Sigma^{m}$
\begin{equation} \label{A_n,theta^+}
    \left\{ \begin{array}{lll}
    \mathrm{D}A_{n}^{+}=0    &  & \mathrm{in} \ \Omega_{0}^{m}, \\
    \mathrm{B} A_{n}^{+}=\hat{\alpha} (\delta_{0}^{2})^{n-1} (\partial_{Y_{3}}  \mathfrak{A}_{n-1}-\dfrac{z'}{r}(\xi) \mathfrak{A}_{n-2}) & & \mathrm{on} \ \Sigma^{m}, \\
    A_{n}^{+}= \delta_{n}^{0} \ g_{\theta} & & \mathrm{on} \ \Gamma^{m} \cup \Gamma_{0}^{+},
    \end{array} \right.
\end{equation}
and
\begin{equation}
    \left\{\begin{array}{lll}
    \partial^{2}_{Y_{3}} \mathfrak{A}_{n}-\big(\frac{1}{\hat{\alpha}}\big)^{2} \mathfrak{A}_{n}= \sum\limits_{p=1}^{n} C_{p} \ \mathfrak{A}_{n-p} & & \mathrm{in} \ \Sigma^{m} \times (0,+\infty), \\
    \mathfrak{A}_{n}=\big( \frac{1}{\hat{\alpha} \delta_{0}^{2}}\big)^{n} A_{n}^{+} & & \mathrm{in} \ \Sigma^{m} \times \{ 0\},
    \end{array} \right.
\end{equation}
where $\delta_{0}=\sqrt{\dfrac{2}{\omega \sigma \mu_{0}}}.$ In (\ref{A_n,theta^+}), $\delta_{n}^{0}$ denotes the Kronecker symbol and we assume that $\mathfrak{A}_{-1}=\mathfrak{A}_{-2}=0.$
In the next section, we make explicit the first asymptotics ($A_{0}^{+},$ $\mathfrak{A}_{0}$) and ($A_{1}^{+},$ $\mathfrak{A}_{1}$) by induction.

\subsection{First terms of the asymptotic expansion}
For $n=0,$ we obtain that $A_{0}^{+}$ solves the problem below
\begin{equation} \label{A_theta,0^+}
    \left\{ \begin{array}{lll}
    \mathrm{D}A_{0}^{+}=0     & & \mathrm{in} \ \Omega_{0}^{m}, \\
    \mathrm{B}A_{0}^{+}=0     & & \mathrm{on} \ \Sigma^{m}, \\
    A_{0}^{+}=g_{\theta}     & & \mathrm{on} \ \Gamma^{m} \cup \Gamma_{0}^{+}.
    \end{array}\right.
\end{equation}
Then according to (\ref{A_theta,0^+}), $\mathfrak{A}_{0}$ solves the following ordinary differential equation (ODE) 
\begin{equation} \label{U_theta,0 ODE}
    \left\{ \begin{array}{lll}
    \partial^{2}_{Y_{3}} \mathfrak{A}_{0} - \big( \frac{1}{\hat{\alpha}} \big)^{2} \mathfrak{A}_{0}=0 & & \mathrm{in} \ \Sigma^{m} \times (0,+\infty), \\
    \mathfrak{A}_{0}=A_{0}^{+} & & \mathrm{on} \ \Sigma^{m} \times \{ 0\}.
    \end{array}\right.
\end{equation}
Then the unique solution of (\ref{U_theta,0 ODE}) such that $\mathfrak{A}_{0} \longrightarrow 0$ as $Y_{3} \longrightarrow +\infty$ writes
\begin{equation} \label{U_theta,0}
    \mathfrak{A}_{0}(\xi,Y_{3})=A_{0}^{+} \big(\tau(\xi) \big) \mathrm{e}^{-\frac{Y_{3}}{\hat{\alpha}}}.
\end{equation}
Next, for $n=1,$ we obtain $A_{1}^{+}$ from (\ref{U_theta,0}) that solves 
\begin{equation} \label{A_theta,1^+}
    \left\{ \begin{array}{lll}
    \mathrm{D}A_{1}^{+}=0     & & \mathrm{in} \ \Omega_{0}^{m}, \\
    \mathrm{B}A_{1}^{+}=-A_{0}^{+}     & & \mathrm{on} \ \Sigma^{m}, \\
    A_{1}^{+}=0     & & \mathrm{on} \ \Gamma^{m} \cup \Gamma_{0}^{+}.
    \end{array}\right.
\end{equation}
Then, according to (\ref{A_theta,1^+}), $\mathfrak{A}_{1}$ solves the following ODE
\begin{equation} \label{U_theta,1 ODE}
    \left\{ \begin{array}{lll}
    \partial^{2}_{Y_{3}} \mathfrak{A}_{1} - \big( \frac{1}{\hat{\alpha}} \big)^{2} \mathfrak{A}_{1}=\frac{1}{\hat{\alpha}} (k+\frac{z'}{r})(\xi) A_{0}^{+} & & \mathrm{in} \ \Sigma^{m} \times (0,+\infty), \\
    \mathfrak{A}_{1}=\frac{1}{\hat{\alpha} \delta_{0}^{2}} A_{1}^{+} & & \mathrm{on} \ \Sigma^{m} \times \{ 0\}.
    \end{array}\right.
\end{equation}
Then the unique solution of (\ref{U_theta,1 ODE}) such that $\mathfrak{A}_{1} \longrightarrow 0$ as $Y_{3} \longrightarrow +\infty$ writes
\begin{equation} \label{U_theta,1}
    \mathfrak{A}_{1}(\xi,Y_{3})=\Bigg[\frac{1}{\hat{\alpha} \delta_{0}^{2}} A_{1}^{+} \big(\tau(\xi) \big)-\frac{Y_{3}}{2} (k+\frac{z'}{r})(\xi) A_{0}^{+}\big(\tau(\xi)\big)  \Bigg] \mathrm{e}^{-\frac{Y_{3}}{\hat{\alpha}}}.
\end{equation}

\section{Elements of proofs of the analytical solutions}\label{Appendix B}
In this section, we will perform the same procedure as in \cite{bermudez2007transient} in order to calculate our analytical solutions (\ref{Analytical global solutions}), (\ref{analytical asymptotics}), and (\ref{analytical impedance solution}). We recall that the electric current flows in the $\Vec{e}_{\theta}$ direction and that is uniformly distributed in the $\Vec{e}_{z}$ direction. Then it results the following expressions of the operators $\rm D$ and $\rm B$ (see for instance (\ref{operators in cylindrical coordinates})):
\begin{equation} \label{operators D & B with hypothesis}
    \begin{array}{l}
    \mathrm{D}(r, z; \partial_{r}, \partial_{z})=\partial^{2}_{r}+\dfrac{1}{r} \partial_{r}-\dfrac{1}{r^{2}},      \\
    \mathrm{B}(\xi; \partial_{r}, \partial_{z})=-\partial_{r}-\dfrac{1}{r}.      
    \end{array}
\end{equation}
Noting that, in our cylindrical case, we have $z'(\xi)=1$ and $r'(\xi)=0$.
Equivalently, we get the following expressions of the operators
\begin{equation} \label{operator with hypothesis of z}
    \begin{array}{l}
   \mathrm{D}(r, z; \partial_{r}, \partial_{z})(\cdot)=\partial_{r} \bigg(\dfrac{1}{r} \big(\partial_{r} \ (r \cdot) \big)\bigg),      \\
    \mathrm{B}(\xi; \partial_{r}, \partial_{z})(\cdot)=- \bigg(\dfrac{1}{r} \big(\partial_{r}(r \cdot)\big)\bigg).      
    \end{array}
\end{equation}
In the following, we will calculate by order the analytical value of $A=(A^{+};A^{-}),$ $A_{0}^{+},$ $A_{1}^{+}$ and $A_{1}^{\varepsilon}.$ 
Taking into account that we only have one connected component and according to (\ref{operator with hypothesis of z}), the first and second equations of the problem (\ref{Appendix problem Atheta general}) become
\begin{equation} \label{eqn1 Atheta}
    \begin{array}{llllllllll}
    &  \partial_{r} \bigg(\dfrac{1}{r} \big( \partial_{r}(r A^{+}) \big) \bigg)(r)=0  & & & & & & \mathrm{if} & R_{1}< r < R_{2},
    \end{array}
\end{equation}
\begin{equation} \label{eqn2 Atheta}
    \begin{array}{llll}
    -\partial_{r} \bigg(\dfrac{1}{r} \big( \partial_{r}(r A^{-}) \big) \bigg)(r)+\mathrm{i} \omega \sigma \mu_{r} \mu_{0} A^{-}=0 & & \mathrm{if} & 0 < r < R_{1},
    \end{array}
\end{equation}
with boundary conditions
\begin{equation} \label{boundary condition 1}
    \begin{array}{lll}
    A^{-}(r) \ \mathrm{is} \  bounded \ & \mathrm{as} & r \longrightarrow 0, \\
    \end{array}
\end{equation}
\begin{equation} \label{boundary condition 2}
    \begin{array}{llll}
    A^{+}(r)=\dfrac{k}{r} & & \mathrm{if} & r=R_{2}.
    \end{array}
\end{equation}
Moreover, we have the following transmission conditions
\begin{equation} \label{transmission condition 1}
    \begin{array}{l}
    A^{-}(R_{1})=A^{+}(R_{1}), \\
    \end{array}
\end{equation}
\begin{equation} \label{transmission condition 2}
    \begin{array}{l}
    \dfrac{1}{r} \big( \partial_{r} (r A^{-}) \big) (R_{1})=\mu_{r} \dfrac{1}{r} \big( \partial_{r} (r A^{+}) \big) (R_{1}).
    \end{array}
\end{equation}
Consequently, the model consists of equations (\ref{eqn1 Atheta}), (\ref{eqn2 Atheta}) with boundary conditions (\ref{boundary condition 1}), (\ref{boundary condition 2}) and interface conditions (\ref{transmission condition 1}), (\ref{transmission condition 2}).  
Simple integration of equation (\ref{eqn1 Atheta}) yields to the following form 
\begin{equation}
    A^{+}(r)= \dfrac{r}{2} a + \dfrac{b}{r},
\end{equation}
where \textit{a, b} and \textit{c} are constants. Now, in order to solve equation (\ref{eqn2 Atheta}), we perform the change of variable $\mathrm{x}=r \gamma,$ in $(0, R_{1}),$ where $\gamma=\sqrt{\omega \sigma \mu_{r} \mu_{0}} e^{\mathrm{i}\frac{\pi}{4}}.$
Then, we get 
\begin{equation}  \label{Bessel equation}
\begin{array}{lll}
     \mathrm{x}^{2} \ \partial^{2}_{\mathrm{x}}\tilde{A}^{-}+ \mathrm{x} \ \partial_{\mathrm{x}} \tilde{A}^{-}-(\mathrm{x}^{2}+1) \tilde{A}^{-}=0 & \mathrm{in} & (0, R_{1}),
\end{array}
\end{equation}
where $\tilde{A}^{-}(\mathrm{x})=A^{-}(\dfrac{\mathrm{x}}{\gamma}).$ Equation (\ref{Bessel equation}) is a Bessel equation for which its general equation is given by
\begin{equation}
    \tilde{A}^{-}(\mathrm{x})=c \mathcal{I}_{1}(\mathrm{x})+d \mathcal{K}_{1}(\mathrm{x}),
\end{equation}
where $\mathcal{I}_{1}$ and $\mathcal{K}_{1}$ are the modified Bessel functions of the first and second kind respectively \cite[Chapter 10 - pages 248--250]{olver2010nist}. By using the boundary condition (\ref{boundary condition 1}), we deduce that $d=0.$
On the other hand, the transmission conditions (\ref{transmission condition 1}) and (\ref{transmission condition 2}), and the boundary condition (\ref{boundary condition 2}) imply the system below
\begin{equation}
    \left\{ \begin{array}{l}
    c \mathcal{I}_{1}(\gamma R_{1})=\dfrac{R_{1}}{1} a + \dfrac{b}{R_{1}},      \\
    c \big[ \gamma \mathcal{I}'(\gamma R_{1})+\mathcal{I}_{1}(\gamma R_{1}) \big] =\mu_{r} R_{1} a, \\
    \dfrac{R_{2}}{2} a + \dfrac{b}{R_{2}}=\dfrac{k}{R_{2}}.
    \end{array} \right.
\end{equation}
Thus, we obtain a linear system of order three for constant unknowns \textit{a, b} and \textit{c}. Solving this linear system, we get
\begin{equation*}
   \begin{array}{l}
    a=\dfrac{k}{R_{1}} \dfrac{g_{1}}{g_{2}}, \\
    b=k- a \dfrac{R_{2}^{2}}{2}, \\
    c =\big[  \mathcal{I}_{1}^{-1}(\gamma R_{1}) \dfrac{(R_{1}^{2}-R_{2}^{2})}{2R_{1}}\big]a + \dfrac{k}{R_{1}} \mathcal{I}_{1}^{-1}(\gamma R_{1}).
    \end{array}
\end{equation*}
Noting that $g_{1}$ and $g_{2}$ are constants defined as follows 
\begin{equation*}
    \begin{array}{l}
    g_{1} = \mathcal{I}_{1}(\gamma R_{1})+\gamma R_{1} \mathcal{I}'_{1}(\gamma R_{1}),       \\
    g_{2} = R_{1} \mu_{r} \mathcal{I}_{1}(\gamma R_{1})-g_{1} \dfrac{(R^{2}_{1}-R^{2}_{2})}{2R_{1}}.      
    \end{array}
\end{equation*} 

In a similar manner, we calculate next the analytical values of $A_{0}^{+},$ $A_{1}^{+}$ and $A_{1}^{\varepsilon}.$
We recall that $A_{0}^{+}$ satisfies the following problem 
\begin{equation} \label{problem Atheta0+ appendix}
    \left\{ \begin{array}{lll}
     \partial_{r} \bigg(\dfrac{1}{r} \big(\partial_{r} (r A_{0}^{+}) \big) \bigg)(r)=0,     &  \mathrm{if} & R_{1} < r < R_{2}, \\
     - \dfrac{1}{r} \big( \partial_{r}(r A_{0}^{+}) \big)(r)=0 & \mathrm{if} & r=R_{1} \\
      A_{0}^{+}= \dfrac{k}{r} & \mathrm{if} & r=R_{2}.
    \end{array}
    \right.
\end{equation}
By a simple integration of the first equation of (\ref{problem Atheta0+ appendix}), we get the following form
\begin{equation}
    A_{0}^{+}(r)=\dfrac{r}{2} a_{0} + \dfrac{b_{0}}{r},
\end{equation}
where $a_{0}$ and $b_{0}$ are constants that are deduced from the boundary conditions and satisfies the linear system of order two below
\begin{equation} \label{a0, b0 system}
    \left\{ \begin{array}{l}
    -a_{0}=0,\\
    \dfrac{R_{2}}{2} a_{0} + \dfrac{b_{0}}{R_{2}}=\dfrac{k}{R_{2}}.
    \end{array} \right.
\end{equation}
Solving the system (\ref{a0, b0 system}), we get
\begin{equation*}
   \begin{array}{lll}
   a_{0}=0 & \mathrm{and} & b_{0}=k.
   \end{array}
\end{equation*}
Hence we deduce that 
\begin{equation} \label{value of A0+}
    A_{0}^{+}(r)=\dfrac{k}{r}.
\end{equation}
Similarly, according to (\ref{A1+ problem}) and the expression of $A_{0}^{+}$ in (\ref{value of A0+}), $A_{1}^{+}$ solves the following system
\begin{equation} \label{problem Atheta1+ appendix}
    \left\{ \begin{array}{lll}
     \partial_{r} \bigg(\dfrac{1}{r} \big(\partial_{r} (r A_{1}^{+}) \big) \bigg)(r)=0     &  \mathrm{if} & R_{1} < r < R_{2}, \\
      -\dfrac{1}{r} \big(\partial_{r} (r A_{1}^{+}) \big)(r)=-\dfrac{k}{r} &  \mathrm{if} & r=R_{1},\\
      A_{1}^{+}(r)= 0 & \mathrm{if} & r=R_{2}.
    \end{array}
    \right.
\end{equation}
By a simple integration of the first equation of the above system (\ref{problem Atheta1+ appendix}), $A_{1}^{+}$ satisfies 
\begin{equation}
    A_{1}^{+}(r)=\dfrac{r}{2} a_{1}+\dfrac{b_{1}}{r}
\end{equation}
where $a_{1}$ and $b_{1}$ solves the following linear system deduced from the boundary conditions in (\ref{problem Atheta1+ appendix}):
\begin{equation}
    \begin{array}{l}
    -a_{1} = -\dfrac{k}{R_{1}},       \\
    \dfrac{R_{2}}{2} a_{1} + \dfrac{b_{1}}{R_{2}} = 0.      
    \end{array}
\end{equation}
Solving this linear system, we get the constants $a_{1}$ and $b_{1}$ as follows
\begin{equation*}
    \begin{array}{lll}
        a_{1}=\dfrac{k}{R_{1}} & \mathrm{and} &   b_{1}=-\dfrac{k R_{2}^{2}}{2 R_{1}}.
    \end{array}
\end{equation*}
As a consequence, the second order asymptotic solution in its analytical form writes
\begin{equation}
    A_{0}^{+}(r)+\dfrac{\varepsilon}{\hat{\alpha}} A_{1}^{+}=\dfrac{k}{r}+\dfrac{\varepsilon}{\hat{\alpha}}(\dfrac{k}{R_{1}} \dfrac{r}{2}-\dfrac{ k R_{2}^{2}}{2 R_{1} } \dfrac{1}{r}).
\end{equation}
Finally, the impedance solution $A_{1}^{\varepsilon}$ satisfying (\ref{impedance model}) is calculated in a similar way as the first asymptotics $A_{0}^{+}$ and $A_{1}^{+}.$

\section*{Acknowledgements}
Authors express their gratitude to Ronan PERRUSSEL for his useful remarks and discussions concerning the numerical part.

\bibliography{bibliography} 
\bibliographystyle{abbrv}

\vspace{5mm}
\noindent Dima ABOU EL NASSER EL YAFI, Victor PÉRON\\
Laboratoire de mathématiques appliquées de Pau, E2S UPPA, CNRS \\
Université de Pau et des pays de l'Adour \\
64000 Pau \\
France \\
E-mails: dima.abou-el-nasser-el-yafi@univ-pau.fr; victor.peron@univ-pau.fr \\

\end{document}